\documentclass[a4paper,10pt]{article}
\usepackage[english]{babel}
\usepackage[utf8]{inputenc}
\usepackage[T1]{fontenc}
\usepackage{textcomp} 
\usepackage{amsmath}
\usepackage{amsfonts}
\usepackage{amsthm}
\usepackage{mathrsfs}
\usepackage{amssymb}
\usepackage{csquotes}
\usepackage{graphicx}
\usepackage{bbm}
\usepackage[mathscr]{euscript}
\DeclareSymbolFont{rsfs}{U}{rsfs}{m}{n}
\DeclareSymbolFontAlphabet{\mathscrsfs}{rsfs}
\usepackage{multirow}
\usepackage{multicol,dsfont}
\usepackage{tabulary}
\usepackage{enumerate}

\theoremstyle{definition}
\newtheorem{defi}{Definition}[section]

\newtheorem{remark}[defi]{Remark}

\theoremstyle{plain}
\newtheorem{proposition}[defi]{Proposition}
\newtheorem{theorem}[defi]{Theorem}
\newtheorem{lemma}[defi]{Lemma}

\NewDocumentCommand\iit{O{\R}O{}mm}{\int_{#1}^{#2} #3 \,\mathrm{d}#4\,}
\DeclareMathOperator{\supp}{supp}

\DeclareMathOperator{\N}{\mathbb{N}}

\DeclareMathOperator{\Z}{\mathbb{Z}}
\DeclareMathOperator{\R}{\mathbb{R}}
\DeclareMathOperator{\E}{\mathbb{E}}
\DeclareMathOperator{\p}{\mathbb{P}}

\newcommand{\norm}[1]{\lVert #1 \rVert}
\newcommand{\floor}[1]{\lfloor #1 \rfloor}

\title{A Wavelet-Type Decomposition of Gaussian Fields with Stationary Increments and Binary Self-Similarity}

\author{Thomas Devos, C\'eline Esser and Laurent Loosveldt\footnote{Corresponding author at l.loosveldt@uliege.be}\\[2mm]
\small Department of Mathematics, University of Liège,\\
\small Allée de la Découverte 12, 4000 Liège, Belgium}
\date{}

\begin{document}
\maketitle

\begin{abstract}
We consider Gaussian random fields with stationary increments whose spectral kernel satisfies a binary self-similarity condition. Using a Lemarié-Meyer wavelet basis, we characterize when such fields admit a wavelet-type expansion with almost surely uniform convergence on compact sets. We then introduce nearly orthogonal test functions which allow us to recover localized combinations of the Gaussian coefficients, up to a negligible error. This enables us to prove the existence of slow, ordinary and rapid points. Finally, we obtain a moving-average representation of the field from the same wavelet decomposition.
\end{abstract}


\section{Introduction and statement of the main results}

In recent years, wavelet-type expansions of stochastic processes have become a crucial tool in the analysis of stochastic systems. 
Wavelet expansions provide sharp regularity estimates \cite{ayacheroueffxiao,Ayachexiao2026,Daw2022,Esser2024,Esser2022} and a natural framework for numerical simulation \cite{ABRY20062326,pipiras05,ayachehamonier,ayachehamonierloosveldt2,abrysellan96}. One of their main advantages comes from the fact that they separate the random and deterministic features of the model: the randomness is encoded in a family of coefficients, while the regularity and localization properties are carried by deterministic wavelet-type functions.

For Gaussian processes and fields, the strategies used to obtain such expansions essentially follow the same procedure. Starting from a stochastic integral representation of the process or field, one expands the deterministic integrand in a suitable wavelet basis, and then transfers this expansion to the stochastic integral through the Wiener isometry. This first step usually gives a random series which converges in $L^2(\Omega)$ for every fixed value of the parameter. However, this convergence is generally too weak for pathwise questions. In order to study pointwise regularity, local oscillations, or sample path properties, one needs almost sure uniform convergence on compact sets.

In most classical examples, this strengthening relies on arguments adapted to the process under consideration. For instance, for the usual fractional Brownian motion, the deterministic functions appearing in the wavelet expansion can be expressed in terms of fractional primitives of the analyzing wavelets \cite{meyerselaantaqqu}. If the wavelet basis is smooth enough, these fractional primitives enjoy fast decay properties, which are the key point in the proof of almost sure uniform convergence. Similar strategies have been developed for fractional Brownian sheets \cite{ayachelegerpontier}, multifractional processes \cite{ayachetaqqu2005,ayache2002}, Hermite processes \cite{ayachehamonierloosveldt}, and stable fields \cite{ayacheshiehxiao,ayacheroueffxiao,ayachelouckx}. Although the technical details depend on the model, these approaches share a common starting point: the stochastic expansion separates deterministic wavelet-type kernels from independent coefficients, so that the pathwise analysis relies both on decay and smoothness estimates for the kernels and on bounds for the random coefficients.

The aim of this paper is to investigate this mechanism in a more systematic way for Gaussian fields with stationary increments. We focus on harmonizable Gaussian fields, whose covariance structure is encoded by a spectral kernel $K$. This framework is general enough to include several classical fractional and anisotropic Gaussian fields, but it also raises a specific difficulty: the kernel $K$ may be singular at the origin. The Lemarié-Meyer basis is particularly well suited to this setting, since the Fourier supports of its elements avoid the origin. The dyadic structure of the wavelet basis naturally leads us to impose a dyadic scaling condition on $K$. Instead of requiring the usual self-similarity with respect to every positive dilation, we assume self-similarity only with respect to powers of two. This binary self-similarity is weaker than the usual homogeneity condition, but it is exactly adapted to the scales appearing in the wavelet expansion. 

We now make this framework precise. Let us fix a probability space $(\Omega,\mathcal F,\mathbb P)$, $n\in\mathbb N$, and let $W$ be a Brownian measure on $\mathbb R^n$, that is, a centered Gaussian family
\[
\{W(A): A\subseteq \mathbb R^n \text{ Borel with } \lambda_n(A)<\infty\}
\]
such that, for any Borel sets $A,B\subseteq\mathbb R^n$ of finite Lebesgue measure,
\begin{equation}\label{eqn:presisometrie}
    \mathbb E\left[W(A)\overline{W(B)}\right]=\lambda_n(A\cap B),
\end{equation}
where \(\lambda_n\) denotes the Lebesgue measure on \(\mathbb R^n\). 
This last relation is the key feature to define the so-called Wiener-Itô integral corresponding to $W$. If $a_1,\dots,a_k$ are complex numbers and $A_1,\dots,A_k$ disjoint Borel sets of finite Lebesgue measure, the integral of the simple function
\[ f := \sum_{j=1}^k a_j \mathds{1}_{A_j}\]
is given by
\[ \int_{\R^n} f(x) \, dW(x) = \sum_{j=1}^k a_j W(A_j).\]
Then, from \eqref{eqn:presisometrie}, we deduce that for any such simple function $f,g$, the identity
\begin{equation}\label{eqn:isometry}
    \E\left[\int_{\R^n} f(x) \, dW(x) \overline{\int_{\R^n} g(x) \, dW(x)}\right]= \langle f,g \rangle_{L^2(\R^n)}
\end{equation}
holds. Since simple functions are dense in $L^2(\R^n)$, \eqref{eqn:isometry} allows to extend the stochastic integral as an isometry between the two Hilbert spaces $L^2(\R^n)$ and $L^2(\Omega)$. Now, we consider a stochastic process of the form 
\[\left\{ X_t := \int_{\R^n} f_t(x) \, dW(x) \right\} _{t \in\R^n}, \]
with, for any $t \in \R^n$, $f_t \in L^2(\R^n)$. Then, if $\{\psi_\lambda\}_\lambda$ is a wavelet basis in $L^2(\R^n)$, for any $t \in \R^n$, one can expand $f_t = \sum_\lambda C_\lambda (t) \,  \psi_\lambda$ where
\begin{equation}\label{eqn:coeff_deterministe}
C_\lambda (t)= \int_{\R^n} f_t(x) \overline{\psi_\lambda(x)} \,dx
 \end{equation}
 and deduce from the isometry \eqref{eqn:isometry}
\begin{equation}\label{eqn:convergencedebase}
    X_t =  \int_{\R^n} f_t(x) \, dW(x)= \sum_\lambda C_\lambda (t)  \int_{\R^n} \psi_\lambda(x) \, dW(x).
\end{equation}
with convergence in $L^2(\Omega)$.  However, as already mentioned, this convergence is generally too weak since it holds only for each \textbf{fixed} $t$ in the random space $L^2(\Omega)$, while, for applications such as regularity analysis or simulations, almost sure uniform convergence on compact sets is much more desirable. Thus, a central issue in the study of wavelet-type expansions of Gaussian processes and fields is to strengthen the weak convergence of the random series \eqref{eqn:convergencedebase}, obtained ``for free'' from the isometry \eqref{eqn:isometry} into the stronger almost sure convergence. To achieve this goal, we remark that the isometry \eqref{eqn:isometry}, together with the orthonormality of the wavelet family in  $L^2(\mathbb R^n)$, also entails that $(\int_{\R^n}\psi_\lambda(x) dW(x))_\lambda$ is a family of independent $\mathcal{N}(0,1)$ random variables.  The problem is therefore reduced to estimating the regularity of the deterministic coefficients $t \mapsto C_\lambda (t) $, 
and combining these estimates with probabilistic bounds on the Gaussian coefficients.

For instance, the usual fractional Brownian motion $\{B^H_t\}_{t \geq 0}$ of Hurst parameter $H \in (0,1)$ can be defined by taking $n=1$ and considering, for any $t \geq 0$, the square integrable function\footnote{We denote by $\Gamma$ the usual Gamma function and, for any $\alpha \in \R$, $(\cdot)^\alpha_+$ the function such that, for any $x \in \R$, $(x)^\alpha_+=x^\alpha$ if $x>0$ and $(x)^\alpha_+=0$ otherwise.}
\begin{equation}\label{eqn:kernetfbm}
   f_t \, : \,  s \mapsto \frac{1}{\Gamma(H+1/2) }((t-s)_+^{H-1/2}-(-s)^{H-1/2}_+).
\end{equation}
Then, from the expression \eqref{eqn:kernetfbm}, we deduce that the functions $t \mapsto C_\lambda^{(t)}$ in \eqref{eqn:convergencedebase} can be expressed as a difference of fractional primitive of $\psi_\lambda$. If the wavelet basis is smooth enough, these fractional primitives enjoy fast decay properties which are the key feature to prove the almost sure uniform convergence of the random series \eqref{eqn:convergencedebase}. We refer to \cite{meyerselaantaqqu} for the detailed wavelet analysis. 

There are two natural ways to provide multi-parameter extension of the fractional Brownian motion: the so-called Levy fractional Brownian motion and Fractional Brownian sheet \cite{herbin}. Both models can be described by means of $n$-dimensional Wiener-Itô integrals and the main techniques to prove the almost sure uniform convergence of their respective wavelet-type expansions are very similar to the one-dimensional case. Indeed, in both settings, the functions $t \mapsto C_\lambda^{(t)}$ can be retrieved from fractional primitives of the orthonormal basis $\{ \psi_\lambda \}_\lambda$. The exact formulation of this last expression depends on how a $n$-dimensional wavelet basis is built, by means of tensor products, from a one-dimensional basis. We refer to \cite{ayachelegerpontier} for a wavelet analysis of fractional Brownian sheet. Several extensions of the usual fractional Brownian field have been proposed in the literature, in order to propose Gaussian fields whose Hölder exponent has the more general form. We refer e.g to \cite{ayachetaqqu2005,ayache2002} for their wavelet-type expansions.

Similar wavelet strategies have also been developed in non-Gaussian settings, including Hermite processes and stable fields; see, for instance, \cite{ayachehamonierloosveldt,ayacheshiehxiao,ayacheroueffxiao,ayachelouckx}.

Despite the strong resemblance between these different constructions,  the corresponding arguments are usually developed separately for each model. This raises the following natural question: given a Gaussian field admitting a stochastic integral representation, how can one systematically obtain a useful wavelet-type series expansion, together with almost sure uniform convergence, from analytic properties of the deterministic kernel?

The aim of the current paper is to provide an answer to this question in the setting of Gaussian processes with stationary increments. Given $n \in \mathbb{N}$, it is well-known \cite{Xiao2021TenLectures} that such a $n$-dimensional field $X=\{ X_t\}_{t \in \R^n}$  may be represented by a so-called spectral measure. If this measure is absolutely continuous with respect to the Lebesgue measure, then $X$ has the harmonizable representation
\begin{equation}\label{eqn:representation}
    X_t = \dfrac{1}{(2\pi)^{n/2}} \int_{\R^n} (e^{it\cdot \xi}-1) K(\xi) \,d\widehat{W}(\xi) \quad t \in \R^n,
\end{equation}
where the notation $t\cdot \xi$ stands for the canonical scalar product of the vectors $t$ and $\xi$ in $\R^n$ and the function $K$ is such that, for any $t \in \R^n$,
\begin{equation}\label{eqn:danslintegrale}
    \xi \mapsto (e^{it\cdot \xi}-1) K(\xi) 
\end{equation}
is in $L^2(\R^n)$. In particular, $K \in L_\text{loc}^2(\R^n \setminus \{0\})$ and $|K|^2$ is integrable at infinity. 
Throughout the paper, we use the Fourier transform convention
\[
\widehat f(\xi)=\int_{\R^n}e^{-ix\cdot\xi}f(x)\,dx,\qquad f(x)=\frac{1}{(2\pi)^n}\int_{\R^n}e^{ix\cdot\xi}\widehat f(\xi)\,d\xi.
\] 
In \eqref{eqn:representation}, $\widehat{W}$ denotes the Fourier transform of the Brownian measure $W$. It is characterized by the identity
\begin{equation}\label{eqn:intwithfouriermeasure}
     \int_{\R^n} \overline{\widehat{f}(\xi)} d\widehat{W}(\xi) = (2\pi)^n\int_{\R^n} \overline{f(x)} d W(x), \qquad f \in L^2(\R^n).
\end{equation}

The square-integrability condition for the functions \eqref{eqn:danslintegrale} allows $K$ to have a singularity at the origin.  To avoid this difficulty, we work with the real-valued Lemarié-Meyer basis of $L^2(\R^n)$
\[\{\psi_{p,j,k}:= 2^{jn/2} \psi_p(2^jx-k) : p \in \{1,\dots, 2^n-1\}, j\in \Z, k \in \Z^n\},\]
for which the functions $\psi_p$ belong to the Schwartz space $\mathscr{S}$ and, for all $p \in \{1,\dots, 2^n-1\}, j\in \Z, k \in \Z^n$, $\widehat{\psi}_{p,j,k}(\xi) = 2^{-jn/2}e^{-i2^{-j}k\cdot \xi}\widehat{\psi}_p(2^{-j}\xi)$ is compactly supported away from the origin:
\begin{equation}\label{eqn:compactsupp}
    \supp \widehat{\psi}_{p,j,k} \subset \left[-\dfrac{2^{j+3}\pi}{3},\dfrac{2^{j+3}\pi}{3}\right]^n\setminus \left( -\dfrac{2^{j+1}\pi}{3},\dfrac{2^{j+1}\pi}{3}\right)^n,
\end{equation}
see e.g. \cite{lemariemeyer,Ayache2019}. Combining \eqref{eqn:convergencedebase} with \eqref{eqn:intwithfouriermeasure}, one gets, for any $t \in \R^n$, with convergence in $L^2(\Omega)$,
    \begin{equation} 
    X_t = (2\pi)^{-n/2}\sum_{p=1}^{2^n-1}\sum_{j\in \Z} \sum_{k\in \Z^n} \varepsilon_{p,j,k}\left(\int_{\R^n} (e^{it\cdot \xi}-1)K(\xi)\widehat{\psi}_{p,j,k}(\xi) \,d\xi\right)
    \label{eq: decomposition de X-step1}
    \end{equation}
where 
\begin{align*}
    \varepsilon_{p,j,k}
    & : = (2\pi)^{-n}\int_{\R^n} \!\overline{\widehat{\psi}_{p,j,k}(t)} \, d\widehat{W}(t) \\
    & =\int_{\R^n}\overline{\psi_{p,j,k}(x)}\,dW(x) , \quad (p,j,k)\in \{1,\dots,2^n-1\}\times \Z\times \Z^n
\end{align*} 

are independent $\mathcal{N}(0,1)$.

An application of the Borel-Cantelli Lemma entails the existence of an event $\Omega^*$ of probability $1$ and a positive random variable $C^*$ of finite moment of every order such that, for all $\omega \in \Omega^*$ and $(p,j,k)\in \{1,\dots,2^n-1\}\times \Z\times \Z^n$,
\begin{equation}\label{eqn:moustique}
    |\varepsilon_{p,j,k}(\omega)| \leq C^*(\omega) \sqrt{\log(3+|j|+|k|)},
\end{equation}
see e.g.  \cite[Lemma 2]{ayachetaqqu2003}. Thus, to ensure the almost sure uniform convergence of the series \eqref{eq: decomposition de X-step1}, it would be desirable that the fast decay property of $\widehat{\psi}$ could be translated into decay properties in $j$ and $k$ for the function
\[ t \mapsto \int_{\R^n} (e^{it\cdot \xi}-1)K(\xi)\widehat{\psi}_{p,j,k}(\xi) \,d\xi.\]
To achieve this goal, we propose to assume the following self-similarity condition for the function $K$.
\begin{defi}
    Let $H \in (0,1)$, the function $K$ is said to be \textit{binary $H$-self-similar} if for all $j\in \Z$ and all $\xi\in \R^n\setminus \{0\}$ 
\[K(2^j\xi) = 2^{j(-H-n/2)}K(\xi).\]
\end{defi}
When $K$ is binary $H$-self-similar, the wavelet coefficients of the harmonizable representation take a simple form. Indeed, if $t \in \R^n$ and $(p,j,k)\in \{1,\dots,2^n-1\}\times \Z\times \Z^n$,
\begin{align}
\int_{\R^n} (e^{it\cdot \xi}-1)K(\xi)\widehat{\psi}_{p,j,k}(\xi) \,d\xi 
&= 2^{jn/2}\int_{\R^n} e^{-ik\cdot \xi}(e^{i2^{j}t\cdot \xi}-1)K(2^j\xi)\widehat{\psi_p}(\xi)\,d\xi \nonumber\\
& = 2^{-jH}\int_{\R^n} e^{-ik\cdot \xi}(e^{i2^{j}t\cdot \xi}-1)K(\xi)\widehat{\psi_p}(\xi)\,d\xi \label{eqn:motivationpourK}
\end{align}
Thus, defining, for all $p \in \{1,\dots,2^n-1\}$ the function
   \begin{equation}\label{eqn:defofpsip}
        \Psi_p \, : \, t \mapsto \dfrac{1}{(2\pi)^{n/2}}\int_{\R^n} e^{it\cdot\xi}K(\xi) \widehat{\psi}_p(\xi) \, d \xi,
   \end{equation} 
one can rewrite the series \eqref{eq: decomposition de X-step1} in the form
        \begin{equation} 
    X_t = \sum_{p=1}^{2^n-1}\sum_{j\in \Z} \sum_{k\in \Z^n} 2^{-jH}\varepsilon_{p,j,k}\left(\Psi_p(2^jt-k)-\Psi_p(-k)\right).
    \label{eq: decomposition de X}
    \end{equation}
Thus, the pathwise convergence of the random series is governed by the decay and regularity properties of the functions $\Psi_p$.

\begin{remark}
    The binary $H$-self-similarity is slightly more general than the usual $H$-self-similarity. Indeed,  there exist functions that are binary $H$-self-similar but not $H$-self-similar. For instance, the function
    \[f(x) = |x|^{-H-n/2} \cos(\log_2(|x|)2\pi)\]
    is clearly binary $H$ self-similar but not self-similar. This can be seen by solving the equation $f(ax) = a^{-H-n/2}f(x)$ with respect to $a>0$ for any $x\in \R^n \setminus\{0\}$. More generally, one can show that a kernel $K$ satisfy the binary self-similarity if  and only if there exists a function $P : \R \times S^{n-1} \to \mathbb{C}$ which is $1$-periodic with respect to the first variable such that
\[K(\xi) = |\xi|^{-H-n/2} P\left(\log_2(|\xi|),\frac{\xi}{|\xi|}\right).\]
    When $P$ does not depend on the angular variable, the resulting field will be isotropic. Dependence of $P$ on the angular variable may lead to anisotropy.
    The binary self-similarity condition is naturally dictated by the wavelet-analysis setting, which relies on a dyadic grid. Note that if $P$ is constant with respect to the first variable then $K$ possesses the usual $H$-self-similarity property and can be written in the form
\[K(\xi) = |\xi|^{-H-n/2}K(\xi/|\xi) = |\xi|^{-H-n/2}\tau(\xi/|\xi|)\]
where $\tau $ is defined on the unit sphere and called the topothesis.
In particular, such a kernel $K$ is entirely determined by its restriction to the unit sphere. This family of fields has been studied intensely by several authors, see e.g. \cite{Baek2014,Bierme2007,Bierme2015, Bierme2009, Bierme2026,Bonami2003,Davy2023,Escande2024,Polisano2015} and the references therein. The framework considered in the present paper therefore includes the dyadic counterpart of these classical homogeneous models, under the regularity and decay assumptions stated below.
\end{remark}

\begin{remark}
The binary $H$-self-similarity has been called by some authors as $H$-\emph{semi-self-similarity} of order $2$ (see e.g. \cite{Maejima1999}).
\end{remark}

\begin{remark}
    Remark that, in general, the field \eqref{eqn:representation} is complex-valued. However, from \eqref{eqn:intwithfouriermeasure} and basic properties of the Fourier transform, one can easily see that the field \eqref{eqn:representation} is real-valued if and only if, for almost every $\xi \in \R^n$, $K(\xi) = \overline{K(-\xi)}$.
\end{remark}

Now, in order to study the almost sure uniform convergence of the random series \eqref{eq: decomposition de X}, we propose to investigate under which conditions the functions in \eqref{eqn:defofpsip} belong to a function space\footnote{Of course, we have $\mathscr{S}=\bigcap_{m \in \N} \mathscr{S}^m$.}
\[\mathscr{S}^\alpha = \left\{f\in C^{\infty}(\R^n) : \forall \gamma \in \N^n \exists C>0 : \forall t\in \R^n |\partial^\gamma f(t)|\leq C (1+|t|)^{-\alpha} \right\},\]
for some $\alpha>0$.

 In fact, we characterize this property in terms of the wavelet coefficients of $K$.  To this end, we introduce the following sequence space.

 \begin{defi}\label{defofsHalpha}
Let $H\in(0,1)$ and $\alpha>0$. We define $s^{H,\alpha}$ as the set of sequences
$
(c_{p,j,k})_{p\in\{1,\dots,2^n-1\},\,j\in\mathbb Z,\,k\in\mathbb Z^n}
$
such that there exists $C>0$ satisfying
\[
|c_{p,0,k}|\le C(1+|k|)^{-\alpha}
\quad \text{ and }\quad 
c_{p,j,k}=2^{-jH}c_{p,0,k}
\]
for every \(p\in\{1,\dots,2^n-1\}\), \(j\in\mathbb Z\), and \(k\in\mathbb Z^n\).
 \end{defi}

The first main theorem of this paper reads as follows.

\begin{theorem}\label{thm:main1}
Let $H \in (0,1)$ and $\alpha>n$ and $K \in L_\text{loc}^2(\R^n \setminus \{0\})$ be such that $|K|^2$ is integrable at infinity. Then the following two assertions are equivalent:
\begin{enumerate}
\item The following three properties hold simultaneously:
\begin{enumerate}[(i)]
\item $K$ admits a binary $H$-self-similar representative;
\item for all $t\in \R^n$, the function $\xi \mapsto (e^{it\cdot \xi}-1)K(\xi)$ belongs to $L^2(\R^n)$;
\item for all $p \in \{1,\dots,2^n-1\}$, the function $\Psi_p$ defined by \eqref{eqn:defofpsip} belongs to $\mathscr{S}^\alpha$.
\end{enumerate}
\item The sequence $(c_{p,j,k})_{p \in \{1,\dots,2^n-1\}, j \in \Z, k \in \Z^n}$ defined by
\begin{equation}\label{eq:c_pjk}
    c_{p,j,k} := \int_{\R^n}{K}(\xi) \widehat{\psi}_{p,j,k}(\xi)\, d\xi
\end{equation}
belongs to $s^{H,\alpha}$.
\end{enumerate}
\end{theorem}

This theorem leads to the following definition.

\begin{defi}\label{def:(H,alpha)condi}
The Gaussian field \eqref{eqn:representation} is said to \textit{satisfy the $(H,\alpha)$-condition} if the sequence \eqref{eq:c_pjk} belongs to $s^{H,\alpha}$.
\end{defi}

In short, if the Gaussian field \eqref{eqn:representation} satisfies the $(H,\alpha)$-condition, then, for every $t \in \R^n$, the representation \eqref{eq: decomposition de X} holds in $L^2(\Omega)$, and the functions $\Psi_p$ appearing in this decomposition belong to $\mathscr{S}^\alpha$. As a consequence, one deduces the almost sure uniform convergence of this wavelet representation.

\begin{theorem}
\label{thm : conv unif de X}
Assume that the field $X$ given by \eqref{eqn:representation}  satisfies the $(H,\alpha)$-condition, with $H \in (0,1)$ and $\alpha >n$. Then, for all $t,s \in \R^n$, we have 
\begin{equation}\label{eqn:boundofsecondmoment}
    \E[|X_t-X_s|^2]\leq C |s-t|^{2H}
\end{equation}
Moreover, the equality \eqref{eq: decomposition de X} holds almost surely uniformly on every compact subset of $\R^n$.
\end{theorem}

As explained at the beginning of this introduction, wavelet-type expansions of stochastic fields provide a powerful tool to investigate sample path regularity. For instance, they were used in \cite{Esser2022} to show that fractional Brownian motion exhibits slow, ordinary and fast points. This phenomenon was first established for Brownian motion by Kahane \cite{kahane}. As an application of the previous results, we show that every Gaussian field satisfying the condition introduced in Definition \eqref{def:(H,alpha)condi} also exhibits these three types of points in its pointwise regularity.

\begin{theorem}
\label{theoreme principal}
Assume that the field $X$ given by \eqref{eqn:representation} satisfies the $(H,\alpha)$-condition with $\alpha>n$, and assume that its spectral kernel is non-degenerate, in the sense that $K\not\equiv0$. Then, almost surely, for every non-empty open set $I\subset \mathbb R^n$, the following properties hold: 
\begin{itemize}
\item[$\bullet$] almost every $t\in I$ is ordinary, in the sense that 
\begin{equation}\label{eqn:limitordinary}
    0<\limsup_{s\to t} \dfrac{|X_t-X_s|}{|t-s|^H \sqrt{\log	\log |t-s|^{-1}}}<+\infty,
\end{equation}
\item[$\bullet$] there exists $t\in I$ which is fast,  in the sense that
\begin{equation}\label{eqn:limitfast}
    0<\limsup_{s\to t} \dfrac{|X_t-X_s|}{|t-s|^H \sqrt{\log |t-s|^{-1}}}<+\infty,
\end{equation}
\item[$\bullet$] there exists $t\in I$ which is slow, in the sense that
\begin{equation}\label{eqn:limitslow}
    0<\limsup_{s\to t} \dfrac{|X_t-X_s|}{|t-s|^H}<+\infty.
\end{equation}
\end{itemize}
\end{theorem}

At this stage, it is interesting to mention that, in the case where $X$ is the usual Fractional Brownian motion, the proof of the finiteness of the limsups in \eqref{eqn:limitordinary}, \eqref{eqn:limitfast} and \eqref{eqn:limitslow} is obtained in \cite{Esser2022} by writing explicitly the increments of the process in terms of the random series \eqref{eq: decomposition de X}. Then, taking advantage of almost sure bounds for the Gaussian family $(\varepsilon_{p,j,k})_{p,j,k}$ and the fast decay properties of $\Psi_p$, one obtains the desired finiteness of the limits. In particular, Theorem \ref{thm:main1} places us in a similar setting for any field satisfying the $(H,\alpha)$-condition. Thus, the finiteness of the limits \eqref{eqn:limitordinary}, \eqref{eqn:limitfast} and \eqref{eqn:limitslow} in this general setting follows from the arguments of \cite[Propositions 3.2, 3.5,3.8] {Esser2022}, modulo minor adaptations to the multidimensional setting.

By contrast, the positivity of the limsups in
\eqref{eqn:limitordinary}, \eqref{eqn:limitfast} and
\eqref{eqn:limitslow} is more delicate. In the case of fractional Brownian
motion, it was obtained in \cite{Esser2024}, see Remark 3.6 therein, by using
the fact that, in this particular setting, the functions involved in the
wavelet representation form a biorthogonal system; see
\cite[Proposition 5.13]{Ayache2019}. This property makes it possible to
recover each Gaussian coefficient $\varepsilon_{p,j,k}\) appearing in
\eqref{eq: decomposition de X} as a Lebesgue integral involving the sample
paths of the fractional Brownian motion; see \cite[Lemma 2.3]{Esser2024}.

In the general framework of the present paper, there is no reason to expect
such a biorthogonal system to exist. We therefore develop a new strategy to
prove the positivity of the limsups in \eqref{eqn:limitordinary},
\eqref{eqn:limitfast} and \eqref{eqn:limitslow} for general fields satisfying
the \((H,\alpha)\)-condition. The idea is to introduce, for each \(p\), a
function denoted by \(\Psi_p^\perp\), which is nearly orthogonal to the family
of functions \(\Psi_{p'}\). More precisely, for
\(p\in\{1,\dots,2^n-1\}\), \(j\in\mathbb Z\) and \(k\in\mathbb Z^n\), the
quantity
\[
\left|\left\langle X, 2^{jn}\Psi_p^\perp(2^j\cdot-k) \right\rangle \right|
\]
is mainly determined by those Gaussian coefficients
\(\varepsilon_{p',j',k'}\) for which \(j'\in\{j-1,j,j+1\}\) and \(k'\) is
close to \(k\). This almost-local extraction of the random coefficients
provides a substitute for the exact biorthogonality available in the
fractional Brownian motion case. We believe that this approach may also be
useful for further investigations of the pointwise regularity and
irregularity of stochastic fields.

Finally, let us point out that harmonizable representation \eqref{eqn:representation} of stochastic fields with stationary increments is canonically obtained from their spectral measures. However, it is sometimes useful to also provide a representation of such fields in the ``time'' domain, that is through integration with respect to the measure $W$ \cite{Baek2014}. For a field with stationary increments, such a
representation is usually called a moving average representation.

Using the tools developed in the proof of Theorem \ref{thm:main1}, we obtain
a moving average representation for any field satisfying the $(H,\alpha)$-condition, by means of the functions $\Psi_p$ in \eqref{eqn:defofpsip}.

\begin{theorem}\label{thm:moving}
     If the field $X$ satisfies the $(H,\alpha)$-condition with $\alpha>n$ and $H<n/2$, then, for every $t\in\R^n$, it admits the moving average representation
     \[X_t=\int_{\R^n} \big(\varphi(u-t)-\varphi(u)\big) \, dW(u).\]
     where 
     \[\varphi(u) = \sum_{p,j,k} 2^{-jH} \Psi_p(-k) \overline{\psi}_{p,j,k}(u) ,\]
     the convergence being uniform on every compact of $\R^n\setminus\{0\}$.
 \end{theorem}

The rest of the paper is devoted to the proofs of the theorems stated in this introduction.

\section{Proof of Theorem \ref{thm:main1}} \label{sec:proofmain1}

This section is devoted to the proof of Theorem \ref{thm:main1}. The direct implication follows from the binary self-similarity of $K$. The main part of the proof is the converse implication. Starting from the assumption that the wavelet coefficients of $K$ belong to $s^{H,\alpha}$, we first obtain a binary $H$-self-similar representative using the Lemarié-Meyer expansion. We then prove that the corresponding functions $\Psi_p$ belong to $\mathscr S^\alpha$. Finally, we show that $\xi\mapsto (e^{it\cdot\xi}-1)K(\xi)$ belongs to $L^2(\mathbb R^n)$ for every $t \in \R^n$, which completes the proof of the theorem. We end the section with a sufficient condition, formulated directly in terms of the local regularity of $K$, which ensures that the functions $\Psi_p$ have the required decay.

\medskip

Let $H \in (0,1)$, $\alpha>n$, and $K \in L_\text{loc}^2(\R^n \setminus \{0\})$ be such that $|K|^2$ is integrable at infinity. We first prove the direct implication in Theorem \ref{thm:main1}. Assume that $K$  is a binary $H$-self-similar function such that, for all $t\in \R^n$, the function $$\xi \mapsto (e^{it\cdot \xi}-1)K(\xi)$$ belongs to $L^2(\R^n)$ and, for all $p\in \{1,\dots, 2^n-1\}$, the functions $\Psi_p$ defined by \eqref{eqn:defofpsip} belong to $\mathscr{S}^\alpha$.
Then the sequence \eqref{eq:c_pjk} belongs to $s^{H,\alpha}$: Indeed, for any $p \in \{1,\dots, 2^n-1\}$, $j\in \Z$ and $k \in \Z^n$, we have
\begin{align*}
\int_{\mathbb R^n} K(\xi)\widehat{\psi}_{p,j,k}(\xi)\,d\xi
&=2^{-jn/2}\int_{\mathbb R^n}e^{-i2^{-j}k\cdot \xi}K(\xi)\widehat{\psi}_p(2^{-j}\xi)\,d\xi  \\
&=2^{jn/2}\int_{\mathbb R^n}e^{-ik\cdot \eta}K(2^j\eta)\widehat{\psi}_p(\eta)\,d\eta  \\
&=(2\pi)^{n/2}2^{-jH}\Psi_p(-k).
\end{align*}
Since $\Psi_p\in\mathscr{S}^\alpha$, we have
$$
|\Psi_p(-k)|\le C(1+|k|)^{-\alpha},
$$
uniformly in $k\in\mathbb{Z}^n$, and therefore the sequence
\eqref{eq:c_pjk} belongs to $s^{H,\alpha}$.

It remains to prove the reverse implication in Theorem \ref{thm:main1}.
Before doing so, we establish two preliminary lemmas.

\begin{lemma}\label{lemma:conv_serie}
Let $\alpha>n$ and $f\in \mathscr{S}^\alpha$. For every $j\in \Z$, the series
\[
\sum_{k\in \Z^n} f(2^jx-k)
\]
converges absolutely for every $x\in\R^n$. Moreover,
\[
\sup_{j\in\Z,\,x\in\R^n}\sum_{k\in\Z^n}|f(2^jx-k)|<\infty.
\]
\end{lemma}
\begin{proof}
For $y\in \R^n$, let us write $\floor{y} = (\floor{y_1},\dots, \floor{y_n})$.
Using the fact that $f\in \mathscr{S}^\alpha$ and setting $k'= k-\floor{2^jx}$, we get
    \[\sum_{k\in \Z^n} |f(2^jx-k)|\leq \sum_{k\in \Z^n} \dfrac{C}{(1+|2^jx-k|)^\alpha} =\sum_{k\in \Z^n} \dfrac{C}{(1+|2^jx-\floor{2^jx}-k|)^\alpha} .\]
    Using Minkowski inequality, we then get 
    \[|2^jx-\floor{2^jx}-k|\geq |k|-|2^jx-\floor{2^jx}|\geq |k|-\sqrt{n}.\]
    It follows that 
    \[\sum_{k\in \Z^n : |k|>\sqrt{n}}\dfrac{C}{(1+|2^jx-\floor{2^jx}-k|)^\alpha}\leq \sum_{k\in \Z^n : |k|>\sqrt{n}}\dfrac{C}{(1+|k|-\sqrt{n})^\alpha}< +\infty\]
    by comparison with the Riemann series.
    On the other hand, we have 
    \[\sum_{k\in \Z^n : |k|\leq \sqrt{n}}\dfrac{C}{(1+|2^jx-\floor{2^jx}-k|)^\alpha}\leq C\#\{k\in\Z^n:|k|\leq\sqrt n\}<\infty,\]
    which concludes the proof.
    \end{proof}

\begin{lemma}
\label{lemma : conv sens l2 loc}
    Let $f \in L^2_{loc} (\R^n\setminus\{0\})$ such that $|f|^2$ is integrable at infinity. Then, we have 
    \[f = (2\pi)^{-n}\sum_{p,j,k} \langle f, \overline{\widehat{\psi}_{p,j,k}}\rangle \overline{\widehat{\psi}_{p,j,k}}\]
    in  $L^2_{loc} (\R^n\setminus\{0\})$.
\end{lemma}

\begin{proof} 
Fix a compact $A\subset \R^n\setminus\{0\}$.
    For every $\ell\in \N$, let $\chi_\ell$ be a $C^\infty$ function that equals $1$ on $\R^n\setminus \left[ -\dfrac{2^{-\ell+1}\pi}{3},\dfrac{2^{-\ell+1}\pi}{3}\right]^n$ and is equal to $0$ on  $\left[ -\dfrac{2^{-\ell}\pi}{3},\dfrac{2^{-\ell}\pi}{3}\right]^n$. 
    It is clear that, for every $\ell\in \N$, $f\chi_\ell \in L^2(\R^n)$ and thus, we have
    \[f\chi_\ell = (2\pi)^{-n}\sum_{p,j,k} \langle f\chi_\ell,\overline{\widehat{\psi}_{p,j,k}}\rangle \overline{\widehat{\psi}_{p,j,k}}\]
    in $L^2(\R^n)$. Now remark that there exists $L>0$ such that, for every $\ell>L$, $\chi_\ell = 1$ on $A$. Fix $\ell>L$, we have
    \[f = f\chi_\ell =  (2\pi)^{-n}\sum_{p,j,k} \langle f\chi_\ell,\overline{\widehat{\psi}_{p,j,k}}\rangle \overline{\widehat{\psi}_{p,j,k}} \]
    in $L^2(A)$. Using the compact support of $\overline{\widehat{\psi}_{p,j,k}}$ and as $0\notin A$, there exists $J>0$ such that for all $j<-J$, $\overline{\widehat{\psi}_{p,j,k}} =0$ on $A$. This allows us to write, in $L^2(A)$,   \[f = f\chi_\ell =  (2\pi)^{-n}\sum_{j>-J}\sum_{p,k} \langle f\chi_\ell,\overline{\widehat{\psi}_{p,j,k}}\rangle \overline{\widehat{\psi}_{p,j,k}}. \]
    By choosing $\ell \geq J$, we have 
    \[\langle f \chi_\ell, \overline{\widehat{\psi}_{p,j,k}}\rangle = \langle f, \overline{\widehat{\psi}_{p,j,k}}\rangle\]
    for all $j>-J$, $p \in \{1,\dots,2^n-1\}$ and $k \in \Z^n$, which concludes the proof.
\end{proof}

\begin{proposition} \label{prop:represK}
Assume that the sequence $(c_{p,j,k} )_{p,j,k}$  defined by \eqref{eq:c_pjk} belongs to $s^{H,\alpha}$. 
Then, $K$ admits the representation 
\[K = (2\pi)^{-n}\sum_{p,j,k} c_{p,j,k} \overline{\widehat{\psi}_{p,j,k}}\]
in $L^2_{loc}(\R^n\setminus\{0\})$. Moreover, the series
converges uniformly on every compact subset of $\R^n\setminus \{0\}$.
In particular, $K$ admits a representation that is binary $H$-self-similar.
\end{proposition}

\begin{proof}
The convergence in $L^2_{loc}(\R^n\setminus \{0\})$ is a straightforward consequence of the previous lemma.
We now prove the locally uniform convergence away from the origin. 
Let $F$ be such a compact. Using the support of the functions $\widehat{\psi}_{p,j,k}$ (see \eqref{eqn:compactsupp}), we already know that the series over $j$ will only have finitely many non-zero terms. Thus it suffices to consider the series over $k\in \Z^n$. Fix $p\in \{1,\dots,2^n-1\}$ and $j\in \Z$ such that $\supp\widehat{\psi}_{p,j,k}\cap F \ne \emptyset.$ Using that the sequence $(c_{p,j,k} )_{p,j,k}$ belongs to $ s^{H,\alpha}$, we have
$$
|c_{p,j,k}| \le C 2^{-jH}(1+|k|)^{-\alpha}.
$$
Moreover,
$$
\sup_{\xi\in F}\left|\widehat{\psi}_{p,j,k}(\xi)\right|\le C_{p,j},
$$
for every $k$. Hence, 
    \[\sup_{\xi\in F}\sum_{k\in \Z^n} |c_{p,j,k} \overline{\widehat{\psi}_{p,j,k}}(\xi)|\leq C_{p,j}2^{-jH}
\sum_{k\in\mathbb Z^n}
\dfrac{C}{(1+|k|)^\alpha}
<+\infty,
\]
since $\alpha>n$. It gives the uniform convergence of the series on $F$.

It remains to prove the binary $H$-self-similarity. Let $\xi \neq 0$ and $j' \in \Z$. Using the representation above and the Definition \ref{defofsHalpha} of $s^{H,\alpha}$,
     \begin{align*}
         K(2^{j'}\xi) & = (2\pi)^{-n}\sum_{p,j,k} c_{p,j,k} 2^{-j\frac{n}{2}}e^{i2^{j'-j}k\cdot \xi}\overline{\widehat{\psi}_p(2^{j'-j}\xi)} \\
         & = (2\pi)^{-n}2^{-j'\frac{n}{2}} \sum_{p,j,k} c_{p,j+j',k} 2^{-j \frac{n}{2}}e^{i2^{-j}k\cdot \xi}\overline{\widehat{\psi}_p(2^{-j}\xi)} \\
         &= (2\pi)^{-n}2^{j'(-H-\frac{n}{2})} \sum_{p,j,k} c_{p,j,k} 2^{-j \frac{n}{2}}e^{i2^{-j}k\cdot \xi}\overline{\widehat{\psi}_p(2^{-j}\xi)} \\
         &=  2^{j'(-H-\frac{n}{2})} K(\xi).
     \end{align*}
\end{proof}

In the sequel, we write $\mathscr{D}_0$ for the set of test functions whose support is away from zero.

\begin{lemma}\label{lemma:sortserie}
Assume that the sequence $(c_{p,j,k} )_{p,j,k}$  defined by \eqref{eq:c_pjk} belongs to $s^{H,\alpha}$.    For any $\varphi \in \mathscr{D}_0$, one has
    \[\int_{\R^n} K(\xi) \varphi(\xi) \, d\xi= (2\pi)^{-n}\sum_{p,j,k} c_{p,j,k} \int_{\R^n} \overline{\widehat{\psi}}_{p,j,k}(x)\varphi(x)\, dx.  \]
\end{lemma}

\begin{proof}
Let $\varphi\in\mathscr D_0$ and set $F=\operatorname{supp}\varphi$. Since
$F$ is a compact subset of $\mathbb R^n\setminus\{0\}$, Proposition
\ref{prop:represK} gives the uniform convergence on  $F$ of the series
\[
(2\pi)^{-n}\sum_{p,j,k}c_{p,j,k}\overline{\widehat{\psi}_{p,j,k}(\xi)}
\]
towards $K(\xi)$. Therefore,
\[
(2\pi)^{-n}\sum_{p,j,k} c_{p,j,k}\overline{\widehat{\psi}_{p,j,k}(\xi)}\varphi(\xi)
\]
also converges uniformly on $F$, and vanishes outside $F$. We may thus exchange the sum and the integral, which gives
\[
\int_{\mathbb R^n}K(\xi)\varphi(\xi)\,d\xi
= (2\pi)^{-n}\sum_{p,j,k} c_{p,j,k}\int_{\mathbb R^n}
\overline{\widehat{\psi}_{p,j,k}(\xi)}\varphi(\xi)\,d\xi.
\]
\end{proof}

\begin{proposition}\label{PsipinSalpha}
Assume that the sequence $(c_{p,j,k} )_{p,j,k}$  defined by \eqref{eq:c_pjk} belongs to $s^{H,\alpha}$. Then, for every $p\in \{1,\dots, 2^n-1\}$, the function $\Psi_p$ defined by \eqref{eqn:defofpsip}
belongs to the space $\mathscr{S}^\alpha$.
\end{proposition}

\begin{proof} 
We apply Lemma \ref{lemma:sortserie} to the test function
$
\eta\mapsto e^{it\cdot\eta}\widehat\psi_p(\eta),
$
whose support is compact and away from the origin. We obtain
\begin{align*}
\Psi_p(t)
&= (2\pi)^{-3n/2} \sum_{p',j,k} c_{p',j,k} \int_{\mathbb R^n} e^{it\cdot\eta} \overline{\widehat{\psi}}_{p',j,k}(\eta) \widehat\psi_p(\eta)\,d\eta  \\
&=(2\pi)^{-3n/2}\sum_{p',j,k}2^{-jn/2} c_{p',j,k}\int_{\mathbb R^n}e^{i\eta\cdot(t+2^{-j}k)}\overline{\widehat\psi_{p'}}(2^{-j}\eta)\widehat\psi_p(\eta)\,d\eta .
\end{align*}
In the sequel, we write, for any $p \in \{1,\dots,2^n-1\}$ and $j \in \Z$, $\psi_{p,j}= 2^{jn/2}\psi_p(2^{j}\bullet)$.  With the change of variable \(\eta=2^j\xi\), the previous equality becomes
 \begin{align}
    \Psi_p(t) 
    &=(2\pi)^{-3n/2}\sum_{p',j,k} 2^{jn/2}c_{p',j,k} \int_{\R^n} e^{i\xi\cdot (2^jt+k)}\overline{\widehat{\psi_{p'}}}(\xi)\widehat{\psi_p}(2^j\xi)\, d \xi \notag\\
     &=(2\pi)^{-3n/2}\sum_{j=-1}^1\sum_{p',k} c_{p',j,k} \int_{\R^n} e^{i\xi \cdot (2^jt+k) }\widehat{\overline{\psi}_{p'}(-\cdot)}\widehat{\psi_{p,-j}}(\xi) \, d\xi\notag\\
     &= (2\pi)^{-n/2}\sum_{j=-1}^1\sum_{p',k} c_{p',j,k}  \overline{\psi}_{p'}( - \cdot)*\psi_{p,-j}(2^jt+k) \label{eq7}
\end{align}
where we use the support properties of the Lemarié-Meyer wavelets \eqref{eqn:compactsupp} that  gives $\widehat{\psi_{p'}}\widehat{\psi_{p,j}}$ vanishes as soon as $|j| >1$.

Concerning the smoothness of $\Psi_p$, we note that, for any $j \in \{-1,0, 1 \}$, $k \in \Z^n$ and $p' \in \{1,\dots,2^n-1\}$, we have, for all $\beta \in \N^n$,
\begin{align*}
    \partial^\beta(\overline{\psi}_{p'}(-\cdot)*\psi_{p,-j}(2^jt+k)) = 2^{j|\beta|}\overline{\psi}_{p'}(-\cdot)*\partial^\beta(\psi_{p,-j})(2^jt+k) =\overline{\psi}_{p'}(-\cdot)*\partial^\beta(\psi_{p})_{-j}(2^jt+k) 
\end{align*}
where we note $\partial^\beta(\psi_{p})_{j} := 2^{jn/2}\partial^\beta (\psi_p)(2^j\cdot)$.
Thus, we obtain from \eqref{eq7} that $\Psi_p \in C^{\infty}(\R^n)$, as the series over $k$ is uniformly convergent since the convolution of Schwartz function is also a Schwartz function and Lemma \ref{lemma:conv_serie}.
It remains to prove the decay property. To this end, we write
, for any $t \in \R^n$,
\begin{align}
    (2\pi)^{n/2}|t|^\alpha |\Psi_p(t)| &\leq |t|^\alpha\sum_{j=-1}^1\sum_{p',k} |c_{p',j,k}  \overline{\psi}_{p'}(-\cdot)*\psi_{p,-j}(2^jt+k)|\notag\\
     &= \sum_{j=-1}^1\sum_{p',k}|c_{p',j,k}|  2^{-j\alpha}|2^jt+k-k|^\alpha|\overline{\psi}_{p'}(-\cdot)*\psi_{p,-j}(2^jt+k)|\notag\\
    &\leq \sum_{j=-1}^1\sum_{p',k} 2^\alpha|c_{p',j,k}|  2^{-j\alpha}|2^jt+k|^\alpha|\overline{\psi}_{p'}(-\cdot)*\psi_{p,-j}(2^jt+k)|\notag\\
    &+\sum_{j=-1}^1\sum_{p',k} 2^\alpha|k|^\alpha|c_{p',j,k}|  2^{-j\alpha}|\overline{\psi}_{p'}(-\cdot)*\psi_{p,-j}(2^jt+k)|, \label{eq:sumforproof}
\end{align}
where we used the fact that, for every $a,b,\alpha >0$, 
\begin{equation}\label{eq:exp_somme}
    (a+b)^\alpha \leq 2^\alpha a^\alpha + 2^\alpha b^\alpha.
\end{equation}

Let us bound the first series in \eqref{eq:sumforproof}. As the convolution product of Schwartz function is a Schwartz function and the sequence $c_{p',j,k}\in s^{H,\alpha}$, we have
\begin{align*}
    \sum_{j=-1}^1\sum_{p',k} 2^\alpha|c_{p',j,k}|  2^{-j\alpha}|2^jt+k|^\alpha|\overline{\psi}_{p'}(-\cdot)*\psi_{p,-j}(2^jt+k)| &\leq C\sum_{j=-1}^1\sum_{p',k} 2^\alpha 2^{-j\alpha}\dfrac{C' 2^{-jH}}{(1+|k|)^\alpha}\\
    &\leq C".
\end{align*}

To conclude the proof, it suffices to bound the second series in \eqref{eq:sumforproof}. As the sequence $c_{p',j,k}$ belongs to $s^{H,\alpha}$ and the convolution product of Schwartz function is a Schwartz function, we have
\begin{align*}
    &\sum_{j=-1}^1\sum_{p',k} 2^\alpha|k|^\alpha|c_{p',j,k}|  2^{-j\alpha}|\overline{\psi}_{p'}(-\cdot)*\psi_{p,-j}(2^jt+k)|\\
    & \leq \sum_{j=-1}^1\sum_{p',k} 2^\alpha|k|^\alpha\dfrac{C2^{-jH}}{(1+|k|)^\alpha}  2^{-j\alpha}|\overline{\psi}_{p'}(-\cdot)*\psi_{p,-j}(2^jt+k)|\\
    &\leq C' \sum_{j=-1}^1\sum_{p',k}2^{-j(\alpha+H)}|\overline{\psi}_{p'}(-\cdot)*\psi_{p,-j}(2^jt+k)|< +\infty.
\end{align*}

The same argument applies to the derivatives of $\Psi_p$, using the fact that 
$
\partial^\beta\bigl(\overline{\psi}_{p'}(-\cdot)*\psi_{p,-j}\bigr)
$
is still in the Schwartz class, uniformly with respect to $j\in\{-1,0,1\}$.
\end{proof}

This last proposition shows that, if $K \in L_\text{loc}^2(\R^n \setminus \{0\})$ is such that $|K|^2$ is integrable at infinity and the sequence \eqref{eq:c_pjk} belongs to $s^{H,\alpha}$, then the functions $\Psi_p$ defined in  \eqref{eqn:defofpsip} belong to $\mathscr{S}^\alpha$. Moreover, Proposition \ref{prop:represK} shows that  $K$ admits a binary $H$-self-similar representative. Thus, to complete the proof of Theorem \ref{thm:main1}, it only remains to prove that, for every $t \in \R^n$, the function $\xi \mapsto (e^{it \cdot \xi}-1) K(\xi)$ belongs to $L^2(\R^n)$.

\begin{proposition}
\label{prop wavelet development via dist}
Assume that the sequence $(c_{p,j,k} )_{p,j,k}$  defined by \eqref{eq:c_pjk} belongs to $s^{H,\alpha}$. Then, for every $t \in \R^n$, the 
series
$$
(2\pi)^{-\frac n2} \sum_{p,j,k} 2^{-jH} \left(\Psi_p(2^jt-k)-\Psi_p(-k)\right)
\overline{\widehat{\psi}_{p,j,k}}
$$
converges in $L^2(\R^n)$. Moreover, for every $t \in \R^n$,
\[
(e^{it\cdot\bullet}-1) K(\bullet) =(2\pi)^{-\frac n2}
\sum_{p,j,k}
2^{-jH}
\left(
\Psi_p(2^jt-k)-\Psi_p(-k)
\right)
\overline{\widehat{\psi}_{p,j,k}}
\]
in $L^2(\R^n)$.
\end{proposition}

\begin{proof}
We first observe that, for any $t \in \R^n$, $j \in \Z$, $k \in \Z^n$ and $p \in \{1,\dots,2^n-1\}$, using \eqref{eqn:motivationpourK} and \eqref{eqn:defofpsip}, we have
\[ (2\pi)^{\frac{n}{2}} 2^{-jH}\left( \Psi_p(2^jt-k) - \Psi_p(-k)\right) = \int_{\R^n} (e^{it \cdot \xi}-1)K(\xi) \widehat{\psi}_{p,j,k}(\xi) \, d \xi. \]
Then, we deduce from Proposition \ref{prop:represK}
\begin{align*}
    (2\pi)^{\frac{n}{2}} 2^{-jH}\left( \Psi_p(2^jt-k) - \Psi_p(-k)\right)= (2\pi)^{-n}\int_{\R^n}\sum_{p',j',k'} (e^{it\cdot \xi}-1) c_{p',j',k'} \overline{\widehat{\psi}_{p',j',k'}}(\xi)\widehat{\psi}_{p,j,k}(\xi) \, d\xi
\end{align*}

 For any $t$ and $\xi \in \R^n$, let us set 
$g(t,\xi) : = (2\pi)^{-n}\sum_{p',j',k'} (e^{it\cdot\xi}-1) c_{p',j',k'}\overline{\widehat{\psi}_{p',j',k'}}(\xi)$.
We first prove that the series $g(t,\bullet)$ defines an element of $L^2(\mathbb R^n)$.

By Plancherel's theorem and the Fourier translation rule, it is enough to prove
that
\[
\sum_{p,j,k}c_{p,j,k}
\bigl(\overline{\psi_{p,j,k}}(\cdot-t)-\overline{\psi_{p,j,k}}\bigr)
\]
converges in $L^2(\R^n)$. Indeed, the inverse Fourier transform of
$(e^{it\cdot\xi}-1)\overline{\widehat\psi_{p,j,k}(\xi)}$ is
$\overline{\psi_{p,j,k}}(-\cdot-t)-\overline{\psi_{p,j,k}}(-\cdot)$,
and reflection preserves the $L^2$ norm. 

Let us fix $p$ and $t \in \R^n$. We aim to show that the series
\begin{equation}\label{eqn:seriesdenormeL2}
    \sum_{j,k} |c_{p,j,k}| \norm{\psi_{p,j,k}(\cdot-t)-\psi_{p,j,k}}_{L^2(\R^n)}
\end{equation}
is convergent. Let us first consider the sub-series of \eqref{eqn:seriesdenormeL2} consisting of the terms with indices $j \geq 0$ and recall that we impose $\alpha>n$. We have, since $(c_{p,j,k})_{p,j,k} \in s^{H,\alpha},$
\begin{align*}
\sum_{j\geq 0}\sum_{k\in\mathbb Z^n} |c_{p,j,k}| & \norm{\psi_{p,j,k}(\cdot-t)-\psi_{p,j,k}}_{L^2(\R^n)} \\ &\leq \sum_{j\geq 0}\sum_{k\in\mathbb Z^n} |c_{p,j,k}| \left(\norm{\psi_{p,j,k}(\cdot-t)}_{L^2(\R^n)}+\norm{\psi_{p,j,k}}_{L^2(\R^n)}\right) \\
& \leq C \sum_{j\geq 0} 2^{-j H}\sum_{k\in\mathbb Z^n} (1+|k|)^{-\alpha} <+\infty.
\end{align*}

For $j<0$, the fundamental theorem of calculus and the translation invariance
of the $L^2$ norm give
\[
\begin{aligned}
\|\psi_{p,j,k}(\cdot-t)-\psi_{p,j,k}\|_{L^2}
&\leq |t|\int_0^1\|\nabla\psi_{p,j,k}(\cdot-\theta t)\|_{L^2}\,d\theta\\
&=|t|\,2^j\|\nabla\psi_p\|_{L^2}.
\end{aligned}
\]
Consequently,
\[
\begin{aligned}
&\sum_{j<0}\sum_{k\in\Z^n}|c_{p,j,k}|
\|\psi_{p,j,k}(\cdot-t)-\psi_{p,j,k}\|_{L^2}\\
&\qquad\leq C|t|\sum_{j<0}2^{j(1-H)}
\sum_{k\in\Z^n}(1+|k|)^{-\alpha}<\infty,
\end{aligned}
\]
since $H<1$ and $\alpha>n$.

\end{proof}

We now have all the ingredients needed to prove Theorem \ref{thm:main1}.

\begin{proof}[Proof of Theorem \ref{thm:main1}]
    Let $H \in (0,1)$ and $\alpha>n$ and $K \in L_\text{loc}^2(\R^n \setminus \{0\})$ be such that $|K|^2$ is integrable at infinity.
    
    Assume first that  $K$  is binary $H$-self-similar,  that, for every $t\in \R^n$, the function $\xi \mapsto (e^{it\cdot \xi}-1)K(\xi)$ belongs to $L^2(\R^n)$ and that, for every $p\in \{1,\dots, 2^n-1\}$, the function $\Psi_p$ defined by \eqref{eqn:defofpsip} belongs to $\mathscr{S}^\alpha$.
    Then, for every $p \in \{1,\dots, 2^n-1\}$, $j\in \Z$ and $k \in \Z^n$, we have already observed that
    \begin{align*}
\int_{\R^n}{K}(x) \widehat{\psi}_{p,j,k}(x)\, dx  =(2\pi)^{n/2}2^{-jH} \Psi_p(-k)
\end{align*}
Since $\Psi_p\in \mathscr{S}^\alpha$, we obtain that the sequence defined by \eqref{eq:c_pjk} is in $s^{H,\alpha}$.

Conversely, assume that the sequence defined by \eqref{eq:c_pjk} is in $s^{H,\alpha}$. By Proposition \ref{prop:represK}, $K$ admits a binary $H$-self-similar representation. By Proposition \ref{PsipinSalpha}, for every $p\in \{1,\dots, 2^n-1\}$, the function $\Psi_p$ belongs to $\mathscr{S}^\alpha$. Finally, Proposition \ref{prop wavelet development via dist} shows that, for every $t \in \R^n$, the function 
$$
\xi\mapsto (e^{it\cdot\xi}-1)K(\xi)
$$
belongs to  $L^2(\R^n)$. This proves the reverse implication and
therefore completes the proof of Theorem \ref{thm:main1}.
\end{proof}

Theorem \ref{thm:main1} \textit{characterizes} the membership of the functions $\Psi_p$,  defined by \eqref{eqn:defofpsip}, to the function space $\mathscr{S}^\alpha$ by means of the wavelet coefficient of $K$. Nevertheless, in practice, these coefficients may be difficult to compute. For this reason, we also provide a more direct sufficient condition, based on
local regularity assumptions on $K$.

Let $\alpha=m+\rho$, with $m\in\N_0$ and $\rho\in(0,1]$. In the sequel, $C^\alpha$ denotes the space of functions that are $m$ times continuously differentiable and whose derivatives of order $m$ are $\rho-$Hölder. 
We now prove that the regularity of $K$ near the support of $\widehat\psi_p$ implies the decay of $\Psi_p$. 
The following lemma is an adaptation of \cite[Proposition 5.10 ]{Ayache2019}. Its proof relies on the
classical link between Hölder regularity and decay of the Fourier transform.

\begin{lemma}
Let $p\in\{1,\dots,2^n-1\}$ and let  $\alpha=m+ \rho$ be a positive real number with $m\in \N_0$ and $\rho \in (0,1]$. If $K$ belongs to $C^\alpha$ in a neighbourhood of $\operatorname{supp}\, \widehat{\psi_p}$, then 
$\Psi_p \in \mathscr{S}^\alpha$.
In particular, if K is $C^\infty$ in this neighbourhood, then $\Psi_p \in \mathscr{S}$.
\end{lemma}

\begin{proof}
For all $t \in \R^n$, let $\delta \in \{-1,1\}^n$ be such that $\delta_j = 1 $ if $t_j\geq 0$ and $\delta_j = -1$ otherwise. 
As, for any $\xi \in \R^n$, one can write
\[t \cdot \xi = \sum_{j=1}^n \delta_j(3+|t_j|)\xi_j-3\delta\cdot \xi, \]
we obtain that, for every $\gamma \in \N^n$, 
\[\partial^\gamma \Psi_p(t) = \dfrac{1}{(2\pi)^{n/2}}\int_{\R^n} \left(\prod_{j=1}^n e^{i\delta_j(3+|t_j|)\xi_j}\right) F_{\gamma,\delta,p}(\xi)\,d\xi\]
where
\[ F_{\gamma,\delta,p} \, : \, \xi \mapsto e^{-3i\delta\cdot\xi}(i\xi)^\gamma\widehat{\psi_p}(\xi) K(\xi),\]
with the notation, for $\xi \in \R^n$, $(i\xi)^\gamma = \prod_{j=1}^n (i\xi_j)^{\gamma_j}$.

Now, remark that $F_{\gamma,\delta,p}$ is $C^\alpha(\R^n)$ and is compactly supported and choose $k\in \{1,\dots, n\}$ such that $|t_k| = \max_{1\leq\ell\leq n} |t_\ell|$. It then follows from integration by parts that,
\[|\partial^\gamma\Psi_p(t)| \leq \dfrac{1}{(2\pi)^{n/2}} (3+|t_k|)^{-m}\left|\int_{\R^n}  \left(\prod_{j=1}^n e^{i\delta_j(3+|t_j|)\xi_j}\right) \partial_{\xi_k}^{m} F_{\gamma,\delta,p}(\xi)\,d\xi \right|.\]

Fix $y\in \R$ such that $e^{i\delta_k (3+|t_k|) y}\neq 1 $ and let $e_k$ denote the $k$th vector of the canonical basis of $\R^n$. By translation invariance,
\begin{align*}
\int_{\R^n}  \left(\prod_{j=1}^n e^{i\delta_j(3+|t_j|)\xi_j}\right) \partial_{\xi_k}^{m} F_{\gamma,\delta,p}(\xi)\,d\xi
= e^{-i\delta_k(3+|t_k|)y}\int_{\R^n}  \left(\prod_{j=1}^n e^{i\delta_j(3+|t_j|)\xi_j}\right) \partial_{\xi_k}^{m} F_{\gamma,\delta,p}(\xi-y e_k)\,d\xi.
\end{align*}
Consequently,
\[
\int_{\R^n}  \left(\prod_{j=1}^n e^{i\delta_j(3+|t_j|)\xi_j}\right) \partial_{\xi_k}^{m} F_{\gamma,\delta,p}(\xi)\,d\xi
= \int_{\R^n}  \left(\prod_{j=1}^n e^{i\delta_j(3+|t_j|)\xi_j}\right)
\dfrac{\partial_{\xi_k}^{m}F_{\gamma,\delta,p}(\xi)-\partial_{\xi_k}^{m}F_{\gamma,\delta,p}(\xi-y e_k)}{1-e^{i\delta_k(3+|t_k|)y}}\,d\xi.
\]
Now, we have
\begin{align*} 
|\partial^\gamma\Psi_p(t)| &\leq \dfrac{1}{(2\pi)^{n/2}} (3+|t_k|)^{-m}\left|\int_{\R^n}  \left(\prod_{j=1}^n e^{i\delta_j(3+|t_j|)\xi_j}\right) \partial_{\xi_k}^{m} F_{\gamma,\delta,p}(\xi)\,d\xi \right|\\
&=\dfrac{1}{(2\pi)^{n/2}} (3+|t_k|)^{-m}\left| \int_{\R^n}  \left(\prod_{j=1}^n e^{i\delta_j(3+|t_j|)\xi_j}\right) \dfrac{\partial_{\xi_k}^{m}F_{\gamma,\delta,p}(\xi)-\partial_{\xi_k}^{m} F_{\gamma,\delta,p}(\xi-y e_k)}{1-e^{i\delta_k (3+|t_k|)y}}\,d\xi \right|\\
&\leq \dfrac{1}{(2\pi)^{n/2}} (3+|t_k|)^{-m}C|y|^{\rho}\dfrac{1}{|1-e^{i\delta_k (3+|t_k|)y}|}
\end{align*}
where we used the compact support and the regularity property of $F_{\gamma,\delta,p}$. By letting $y=1/(3+|t_k|)$, it follows that
\[|\partial^\gamma\Psi_p(t)|\leq \dfrac{1}{(2\pi)^{n/2}} (3+|t_k|)^{-\alpha}C\dfrac{1}{|1-e^{i\delta_k}|}\leq C'\dfrac{1}{(2\pi)^{n/2}} (3+|t_k|)^{-\alpha}.\]
Finally, as $|t_k|\geq |t_j|$ for every $j\in \{1,\dots, n\}$, we know that $|t|\leq \sqrt{n}|t_k|$ and thus 
\[|\partial^\gamma\Psi_p(t)|\leq C'\dfrac{1}{(2\pi)^{n/2}} (3+|t_k|)^{-\alpha}\leq C'\dfrac{1}{(2\pi)^{n/2}} \left(3+\dfrac{|t|}{\sqrt{n}}\right)^{-\alpha},\]
which suffices to conclude.
\end{proof}

\section{Proof of Theorem \ref{thm : conv unif de X}}

Any Gaussian field $X$ satisfying the $(H,\alpha)$-condition admits the wavelet-type expansion \eqref{eq: decomposition de X}, where the functions $\Psi_p$ have the required fast decay. We use this decay to prove the second-moment estimate and the almost sure uniform convergence of the random series \eqref{eq: decomposition de X}.

\begin{proof}[Proof of Theorem \ref{thm : conv unif de X}]

Concerning the first part of the theorem, we have to show that,
\[\E[|X(t)-X(s)|^2]\leq C |t-s|^{2H}.\]
Using the independence of the random variables $(\varepsilon_{p,j,k})_{p,j,k}$, we get
\begin{equation}\label{eq:moment2}
    \E[|X(t)-X(s)|^2]= \sum_{p=1}^{2^n-1}\sum_{j\in \Z}\sum_{k\in \Z^n} 2^{-2jH}\left|\Psi_p(2^jt-k)-\Psi_p(2^js-k)\right|^2.
\end{equation}
Fix $p\in \{1,\dots, 2^n-1\}$, $t,s \in \R^n$ and let $\nu \in \Z$ be such that 
\[2^{-\nu}\leq |t-s|< 2^{1-\nu}.\]

If $j\leq \nu$, we set 
\[
F_{p,j,t}(x) = \sum_{k\in \Z^n}2^{-2Hj} |\Psi_p(2^j t-k)-\Psi_p(2^jx-k)|^2
\]
Using the decay properties of $\Psi_p$, the real-valued function
$F_{p,j,t}$ is continuously differentiable. By the mean value theorem,
there exists $x\in S[s,t]$ such that
\begin{align*}
|F_{p,j,t}(s)-F_{p,j,t}(t)|
&\leq 2^{(1-2H)j+1}|t-s|
\sum_{k\in\Z^n}|\Psi_p(2^jt-k)-\Psi_p(2^jx-k)|
|\nabla\Psi_p(2^jx-k)|\\
&\leq C2^{(1-2H)j}|t-s|
\sum_{k\in\Z^n}|\Psi_p(2^jt-k)-\Psi_p(2^jx-k)|.
\end{align*}
For each $k$, choose $x'_k\in S[t,x]$ where
$|\nabla\Psi_p(2^j\cdot-k)|$ attains its maximum on the segment.
The mean-value inequality then gives
\[
|\Psi_p(2^jt-k)-\Psi_p(2^jx-k)|
\leq 2^j|t-s|\,|\nabla\Psi_p(2^jx'_k-k)|.
\]

Since $x'_k\in S[s,t]$ and $j\leq\nu$, we have $|2^jx'_k-2^jt|\leq 2^j|t-s|\leq 2$. Hence, by the decay of $\nabla\Psi_p$,
\[
|\nabla\Psi_p(2^jx'_k-k)|\leq C(1+|2^jt-k|)^{-\alpha}.
\]
Together with Lemma \ref{lemma:conv_serie}, this gives
\[
|F_{p,j,t}(s)-F_{p,j,t}(t)| \leq C 2^{(2-2H)j}|t-s|^2 \sum_{k \in \Z^n}(1+|2^jt-k|)^{-\alpha} \leq C' 2^{(2-2H)j}|t-s|^2.
\]
Thus
\begin{equation}\label{eq:jleqnu}
    \sum_{j\le \nu}\sum_{k\in\mathbb Z^n} 2^{-2Hj} \left|
\Psi_p(2^jt-k)-\Psi_p(2^js-k) \right|^2
\le C'|t-s|^2 2^{(2-2H)\nu} \leq 
C'|t-s|^{2H}
\end{equation}
by the choice of $\nu$. 

If $j> \nu$, by Lemma \ref{lemma:conv_serie}, we have 
\begin{align*}
\sum_{k\in \Z^n} 2^{-2jH}\left|\Psi_p(2^jt-k)-\Psi_p(2^js-k)\right|^2
&\leq 2^{-2Hj}\sum_{k\in \Z^n}4(|\Psi_p(2^jt-k)|^2+|\Psi_p(2^js-k)|^2)\\
&\leq C 2^{-2jH},
\end{align*}
hence
\begin{equation}\label{eq:jgeqnu}
\sum_{j>\nu}\sum_{k\in \Z^n} 2^{-2jH}\left|\Psi_p(2^jt-k)-\Psi_p(2^js-k)\right|^2 
\leq \sum_{j>\nu} C 2^{-2jH} 
\leq C' 2^{-2\nu H} 
\leq C'' |t-s|^{2H}.
\end{equation}

Finally, putting together \eqref{eq:moment2}, \eqref{eq:jleqnu} and \eqref{eq:jgeqnu}, we obtain the first part of the result. For the uniform convergence of \eqref{eq: decomposition de X}, we apply the Gaussian version of the It\^o--Nisio theorem on every cube, using the convergence of the series in $L^2(\Omega)$ for every fixed $t$, together with the previous estimate.
\end{proof}

Since the increments are Gaussian, \eqref{eqn:boundofsecondmoment} implies that, for every $q\geq2$, \[\E[|X_t-X_s|^q]\leq C_q|t-s|^{qH}.\] Choosing $q$ sufficiently large and applying Kolmogorov's continuity theorem, we deduce that $X$ has a version which is locally $(H-\delta)$-Hölder continuous, for every $\delta\in(0,H)$. In the next section, we show, among other things, that this regularity is optimal.

\section{Proof of Theorem \ref{theoreme principal}}\label{sec:irregularity}

Throughout this section, we fix a Gaussian field $X$ satisfying the
$(H,\alpha)$-condition under the assumptions of Theorem \ref{theoreme principal}. As already explained, it remains only to prove the
positivity of the three limits \eqref{eqn:limitordinary},
\eqref{eqn:limitfast} and \eqref{eqn:limitslow}. 

Let us first observe that all the work resides in analysing the high frequencies parts of the decomposition  \eqref{eq: decomposition de X}. Indeed, for the low frequencies parts, we have, for any multi-index $\beta$ such that $|\beta|=1$,
\[\sum_{j<0} \sum_{k\in \Z^n}2^{-jH}\varepsilon_{p,j,k}\partial^\beta\left(\Psi_p(2^jt-k)-\Psi_p(-k)\right) = \sum_{j<0} \sum_{k\in \Z^n}2^{j(1-H)} \varepsilon_{p,j,k}\partial^\beta(\Psi_p)(2^jt-k)\]
which, thanks to \eqref{eqn:moustique}, is almost surely locally uniformly convergent and thus the low-frequency part of the process is differentiable. Since $H<1$, its increments are $O(|t-s|)$ and are therefore negligible with respect to each of the moduli in \eqref{eqn:modofcont}.

Consequently, it is enough to analyse the high-frequency parts $f_p$, $p\in\{1,\dots,2^n-1\}$,
defined by 
\begin{equation}\label{eqn:defoffp}
    f_p = \sum_{j\geq 0}\sum_{k\in \Z^n}2^{-jH}\varepsilon_{p,j,k}\left(\Psi_p(2^j \bullet-k)-\Psi_p(-k)\right).
\end{equation}
We also set $f:=\sum_{p=1}^{2^n-1}f_p$, so that $X-f$ is differentiable.
In the sequel, we will use the following submultiplicative moduli of continuity (after a harmless modification away from the origin):
\begin{equation}\label{eqn:modofcont}
\sigma_1(x)=x^H\sqrt{\log\log\left(\frac{1}{x}\right)}, 
\qquad
\sigma_2(x)=x^H\sqrt{\log\left(\frac{1}{x}\right)}, 
\qquad
\sigma_3(x)=x^H.
\end{equation}

For every $p \in \{1,\dots,2^n-1\}$, we define
\[
\Psi^{\perp}_p(t)
= (2\pi)^{-n/2} \int_{\R^n} e^{it\cdot \xi}|\xi|^{H+n/2}\widehat{\psi}_p(\xi) \, d\xi.
\]
As before, we also write
\[
\Psi^\perp_{p,j,k}
:=
2^{jn}\Psi_p^\perp(2^j\cdot-k).
\]
For convenience, we also write $\Psi_{p,j,k}:=\Psi_p(2^j\cdot-k)$; thus the functions $\Psi_{p,j,k}$ are not $L^2$-normalized, while the test functions $\Psi_{p,j,k}^{\perp}$ are normalized in $L^1$.
These functions will play the role of test functions adapted to the
decomposition \eqref{eq: decomposition de X}. More precisely, we will show
that they are ``nearly orthogonal'', in a sense made precise below, to the
functions $\Psi_p$ appearing in the expansion of $X$.

First, the next proposition relates the increments of the high-frequency part $f$ to its scalar products with the test functions $\Psi^\perp_p$.
Recall that $k_j(t)$ denotes the unique $k\in\Z^n$ such that $t\in 2^{-j}(k+[0,1)^n)$. We set $\lambda_{j,k}:=2^{-j}(k+[0,1)^n)$ and $\lambda_j(t):=\lambda_{j,k_j(t)}$.

\begin{proposition}
\label{prop : maj produit scalaire}
Let $\sigma$ be a submultiplicative modulus of continuity and let $t\in\R^n$. For every $\gamma>0$, there exists an event of probability $1$ on which there exists a random variable $C^*>0$ such that, for every $j$ sufficiently large and every $p,p'\in\{1,\dots,2^n-1\}$,
\[
|\langle f_{p'},\Psi^\perp_{p,j,k_j(t)}\rangle|
\leq C^*\left(\sup\left\{\frac{|f_{p'}(t)-f_{p'}(s)|}{\sigma(|t-s|)}:0<|t-s|<c2^{-j/2}\right\}\sigma(2^{-j})+2^{-j\gamma}\right).
\]
\end{proposition}

\begin{proof}
The proof follows the same lines as that of \cite[Proposition 2.6]{Esser2024}, replacing the function $\Psi$ by $\Psi_p^\perp$.
\end{proof}

The exceptional event in this argument can be chosen independently of $t$: once the almost sure bound \eqref{eqn:moustique} is fixed, the estimate is deterministic. The random constant may of course depend on $t$. Since the number of components is finite, the same argument applies to $f=\sum_{p'=1}^{2^n-1}f_{p'}$. Using this proposition, our problem now reduces to control, for $t \in \R^n$, the quantities $\langle f, \Psi^{\perp}_{p,j,k_j(t)}\rangle$, for large $j \in \N$. Recall that the $(H,\alpha)$-condition for the field $X$ implies that $\Psi_p \in \mathscr{S}^\alpha$ for every $p\in \{1,\dots, 2^n-1\}$.

Using the expression $\eqref{eqn:defoffp}$ and the fact that $\int_{\R^n}\Psi^\perp_p=0$, we see that it is pertinent to study the behaviour of the quantities

\begin{align*}
\langle \Psi_{p,j,k},\Psi_{p',j',k'}^{\perp}\rangle
&= (2\pi)^{-n}\int_{\R^n}2^{-jn}e^{-i2^{-j}k\cdot \xi}e^{i2^{-j'}k'\cdot \xi}\widehat{\Psi_p}(2^{-j}\xi)\overline{\widehat{\Psi^{\perp}_{p'}}}(2^{-j'}\xi)\, d\xi\\
&= \int_{\R^n} 2^{-jn}e^{-i2^{-j}k\cdot \xi}e^{i2^{-j'}k'\cdot \xi} K(2^{-j}\xi) \widehat{\psi}_{p}(2^{-j}\xi)|2^{-j'}\xi|^{H+n/2} \overline{\widehat{\psi_{p'}}}(2^{-j'}\xi)\,d\xi
\end{align*}
for $p,p' \in \{1,\dots,2^n-1\}$, $j,j' \in \Z$ and $k,k' \in \Z^n$.

If $|j-j'|>1$ then, by the support properties of the wavelets, the scalar product is equal to $0$. If $|j-j'|\leq 1$, the change of variables $\eta = 2^{-j} \xi$ gives
\begin{itemize}
\item[$\bullet$] \underline{if $j=j'$ }
\[\langle \Psi_{p,j,k},\Psi_{p',j,k'}^{\perp}\rangle = (2\pi)^{-n}\int_{\R^n} e^{-i(k-k')\cdot \eta} \widehat{\Psi}_p(\eta) \overline{\widehat{\Psi^{\perp}}_{p'}}(\eta)\, d\eta :=\widehat{\Lambda}_{p,p'}(k-k')\]
\item[$\bullet$] \underline{if $j'=j+1$}
\[\langle \Psi_{p,j,k},\Psi_{p',j+1,k'}^{\perp}\rangle = (2\pi)^{-n}\int_{\R^n} e^{-i(k-k'/2)\cdot \eta} \widehat{\Psi}_p(\eta) \overline{\widehat{\Psi^{\perp}}_{p'}}(\eta/2)\, d\eta :=\widehat{\Lambda^+}_{p,p'}(k-k'/2)\]
\item[$\bullet$] \underline{if $j'=j-1$}
\[\langle \Psi_{p,j,k},\Psi_{p',j-1,k'}^{\perp}\rangle = (2\pi)^{-n}\int_{\R^n} e^{-i(k-2k')\cdot \eta} \widehat{\Psi}_p(\eta) \overline{\widehat{\Psi^{\perp}}_{p'}}(2\eta)\, d\eta :=\widehat{\Lambda^-}_{p,p'}(k-2k').\]
\end{itemize}
 Remark that the three $\widehat{\Lambda}$ functions are well defined as $\Psi$ and $\Psi^{\perp}$ have sufficiently rapid decay. The decay of the lambda functions will be critical in our development. This motivates the following lemma.

\begin{lemma}
    For every $p,p'\in \{1,\dots, 2^n-1\}$, the functions $\widehat{\Lambda}_{p,p'}, \widehat{\Lambda^+}_{p,p'}$ and $\widehat{\Lambda^-}_{p,p'}$ are in $\mathscr{S}^\alpha$.
\end{lemma}

\begin{proof}
    Fix $p,p'\in \{1,\dots, 2^n-1\}$. We only show the result for $\widehat{\Lambda}_{p,p'}$ as the other functions are treated similarly. Using the fact that if $f,g\in L^2(\R^n)$, then $(2\pi)^n f*g(t) = \int_{\R^n} e^{it\cdot \xi} \widehat{f}(\xi) \widehat{g}(\xi)\, d\xi$, we get that 
    \[\widehat{\Lambda}_{p,p'}(t) = \Psi_p*(\overline{\Psi_{p'}^{\perp}}(-\cdot))(-t).\]
    Set $h:=\Psi_p*(\overline{\Psi_{p'}^{\perp}}(-\cdot))$. Since $\widehat\Lambda_{p,p'}(t)=h(-t)$, it is enough to prove that $h\in\mathscr S^\alpha$. Using the fact that $\Psi_p \in \mathscr{S}^\alpha$ and $\Psi_{p'}^\perp\in \mathscr{S}$, we have $h\in C^\infty(\R^n)$. 
    Concerning the decay property, the derivatives are treated similarly. Let $\Psi_{p'}^* = \overline{\Psi_{p'}^\perp}(-\cdot)$. Using \eqref{eq:exp_somme}, the triangular inequality and the decay of $\Psi_p$ and $\Psi_{p'}^*$, we get, for every $t\in \R^n$,
    \begin{align*}
        |t|^\alpha|h(t)|&= |t|^\alpha\left|\int_{\R^n} \Psi_p(t-y)\Psi_{p'}^*(y)\,dy\right|\\
        &= \left|\int_{\R^n}|t-y+y|^\alpha \Psi_p(t-y)\Psi_{p'}^*(y)\,dy\right|\\
        &\leq 2^\alpha \int_{\R^n} |t-y|^\alpha|\Psi_p(t-y)||\Psi_{p'}^*(y)|\,dy+2^\alpha \int_{\R^n}|\Psi_p(t-y)|  |y|^\alpha|\Psi_{p'}^*(y)|\,dy\\
        &\leq C \int_{\R^n} |\Psi_{p'}^*(y)|\, dy + C' \int_{\R^n} |\Psi_{p}(y)|\, dy \\
        &<+\infty.
    \end{align*}
    \end{proof}

\begin{remark}
    This last lemma justifies the term ``nearly orthogonal'' used when we introduced the functions $\Psi_p^{\perp}$. Indeed, the function $\Psi_{p,j,k}$ is orthogonal to ``most'' of the test functions $\Psi^\perp_{p',j',k'}$, as soon as $|j-j'|>1$ and, if not, the value of the scalar product $\langle \Psi_{p,j,k},\Psi_{p',j,k'}^{\perp}\rangle$ rapidly decays when $|k-k'| \to + \infty$.

    Let us also mention that, for the particular normalization $K(\xi)=|\xi|^{-H-n/2}$ of the isotropic fractional Brownian kernel, one actually has, for any  $p,p' \in \{1,\dots, 2^n-1\}, j,j'\in \N$ and $k,k'\in \Z^n$,
    \[\langle \Psi_{p,j,k},\Psi_{p',j',k'}^{\perp}\rangle = (2\pi)^n\,\delta_{(p,j,k),(p',j',k')} \]
    which yields, for any $t \in \R^n$,
    \[ \langle X,\Psi_{p,j,k_j(t)}^{\perp}\rangle = (2\pi)^n\,2^{-jH} \varepsilon_{p,j,k_j(t)} \]
For a spectral kernel $K(\xi)=c_H|\xi|^{-H-n/2}$, both identities above carry an additional factor $c_H$. This fact, combined with Proposition \ref{prop : maj produit scalaire}, is a key feature in the proof of Theorem \ref{theoreme principal} for (multi)fractional Brownian motion \cite{Esser2024,Esser2022}. Unfortunately, in our general setting, there is no reason to always be able to set up such a bi-orthogonal system. It is the reason why we propose to work with the ``nearly orthogonal'' property. However, one can easily show that if $K$ is smooth and $K(t) \neq 0$ for every $t\in \R^n\setminus\{0\}$, then it would be possible to set a bi-orthogonal system by taking
    \[\Psi^{\perp}_p(t) =(2\pi)^{-n/2} \int_{\R^n} e^{it\cdot \xi}\dfrac{1}{(2\pi)^n\overline{K(\xi)}}\widehat{\psi}_p(\xi) \, d\xi. \]
    The conjugate in the denominator is necessary for the sesquilinear $L^2$ scalar product. Then, the consideration of the current section would be greatly simplified, and reduce to the arguments developed in the papers \cite{Esser2022,Esser2024}.
\end{remark}

Since $K\not\equiv0$ and is binary self-similar, the Fourier supports of the Lemari\'e--Meyer wavelets imply that there exists $p\in\{1,\dots,2^n-1\}$ such that $K\widehat\psi_p\not\equiv0$. We fix such an index $p$ for the rest of the proof. Let $j\geq1$ and $k\in \Z^n$. We have, for all $\delta >0$, 

\begin{align*}
    |\langle f, \Psi_{p,j,k}^{\perp}\rangle| \geq 2^{-jH} (|Z_{j,k}(\delta)| - |\widetilde{Z}_{j,k}(\delta)|)
\end{align*}

where we set
\begin{align*}
Z_{j,k}(\delta)
&:=\sum_{p'=1}^{2^n-1}\Bigg(
\sum_{|k'-k|<2^{j\delta}}
\varepsilon_{p',j,k'}\widehat\Lambda_{p',p}(k'-k)
+2^{-H}\sum_{|k'-2k|<2^{j\delta}}
\varepsilon_{p',j+1,k'}\widehat\Lambda^-_{p',p}(k'-2k)\\
&\hspace{3.5cm}+2^H\sum_{|k'-k/2|<2^{j\delta}}
\varepsilon_{p',j-1,k'}\widehat\Lambda^+_{p',p}(k'-k/2)\Bigg)
\end{align*}
and
\begin{align*}
\widetilde Z_{j,k}(\delta)
&:=\sum_{p'=1}^{2^n-1}\Bigg(
\sum_{|k'-k|\geq 2^{j\delta}}
\varepsilon_{p',j,k'}\widehat\Lambda_{p',p}(k'-k)
+2^{-H}\sum_{|k'-2k|\geq 2^{j\delta}}
\varepsilon_{p',j+1,k'}\widehat\Lambda^-_{p',p}(k'-2k)\\
&\hspace{3.5cm}+2^H\sum_{|k'-k/2|\geq 2^{j\delta}}
\varepsilon_{p',j-1,k'}\widehat\Lambda^+_{p',p}(k'-k/2)\Bigg).
\end{align*}
The factors $2^{-H}$ and $2^H$ come from the coefficients
$2^{-(j+1)H}$ and $2^{-(j-1)H}$ in \eqref{eqn:defoffp}.

We first show that the random variable $\widetilde{Z}_{j,k}(\delta)$, $j \in \N$ and $k \in \Z^n$, are negligible (with regard to the moduli of continuity of interest in \eqref{eqn:modofcont}). 

\begin{lemma}
For every compact set $Q\subset\R^n$ and every $\delta>0$, there exists an event of probability $1$ on which there are a random variable $C_Q^*>0$ and $J_Q\in\N$ such that, for every $j\geq J_Q$ and every $k\in\Z^n$ with $\lambda_{j,k}\cap Q\neq\emptyset$,
\[|\widetilde{Z}_{j,k}(\delta)| \leq C_Q^* j 2^{-j\delta(2\alpha-n)/2}.\]
\end{lemma}

\begin{proof}
First remark that $\widetilde{Z}_{j,k}(\delta)$ is a complex-valued Gaussian random variable as the $\hat{\Lambda}$ function are complex-valued functions. As, for any $z\in \mathbb{C}$, $|z| \leq \sqrt{2} \max(|\Re(z)|, |\Im(z)|)$, it suffices to show the result for $\Re\widetilde{Z}_{j,k}(\delta)$ and $\Im\widetilde{Z}_{j,k}(\delta)$. 
We focus on the random variable 
\[\sum_{|k-k'|\geq 2^{j\delta}}  \varepsilon_{p',j,k'} \Re\widehat{\Lambda}_{p',p}(k'-k),\]
as the imaginary part and the two other terms in the definition of $\widetilde{Z}_{j,k}(\delta)$ can be treated similarly. We have, since $\widehat{\Lambda}_{p',p}\in \mathscr{S}^\alpha$,

\begin{align*}
\left\|\sum_{|k-k'|\geq 2^{j\delta}} \varepsilon_{p',j,k'} \Re\widehat{\Lambda}_{p',p}(k'-k)\right\|_{L^2(\Omega)}^2&\leq  \sum_{|k-k'|\geq 2^{j\delta}} \left|\widehat{\Lambda}_{p',p}(k'-k)\right|^2\\
&\leq C^2\sum_{|k|\geq 2^{j\delta}} \dfrac{1}{(1+|k|)^{2\alpha}}\\
&\leq C^2 2^{-j\delta(2\alpha-n)}.
\end{align*}

Now recall that there exists $c>0$ such that for every Gaussian random variable $X$ and $y \geq 2$, $\mathbb{P}(|X| \geq  y \|X\|_{L^2(\Omega)}) \leq \exp(-c y^2)$ (see \cite[Theorem 6.7]{Janson1997}). Thus we get, for any $j \geq 2$,
\begin{align*}
    &\mathbb{P}\left(\left|\sum_{|k-k'|\geq 2^{j\delta}} \varepsilon_{p',j,k'} \Re\widehat{\Lambda}_{p',p}(k'-k)\right| \geq  j C2^{-j\delta(2\alpha-n)/2}\right)\\&\leq \mathbb{P}\left(\left|\sum_{|k-k'|\geq 2^{j\delta}} \varepsilon_{p',j,k'} \Re\widehat{\Lambda}_{p',p}(k'-k)\right| \geq  j \left\|\sum_{|k-k'|\geq 2^{j\delta}} \varepsilon_{p',j,k'} \Re\widehat{\Lambda}_{p',p}(k'-k)\right\|_{L^2(\Omega)}\right) \\ 
    &\leq \exp(-c j^2),
\end{align*}
which implies 
\[\sum_{j\in \N} \mathbb{P}\left(\left|\sum_{|k-k'|\geq 2^{j\delta}} \varepsilon_{p',j,k'} \Re\widehat{\Lambda}_{p',p}(k'-k)\right| \geq  j C2^{-j\delta(2\alpha-n)/2}\right)<+\infty. \]
For a fixed compact $Q$, the number of indices $k$ such that $\lambda_{j,k}\cap Q\neq\emptyset$ is at most $C_Q2^{jn}$. Hence the union bound gives a summable upper bound $C_Q2^{jn}e^{-cj^2}$. Finally, using the Borel-Cantelli Lemma, we deduce the asserted estimate, uniformly over the indices $k$ such that $\lambda_{j,k}\cap Q\neq\emptyset$.
\[
\left|\sum_{|k-k'|\geq 2^{j\delta}} \varepsilon_{p',j,k'} \Re \widehat{\Lambda}_{p',p}(k'-k)\right|
\leq C_Q^* j 2^{-j\delta(2\alpha-n)/2},
\qquad j\geq J_Q,\quad \lambda_{j,k}\cap Q\neq\emptyset.
\]
which concludes the proof.
\end{proof}

Using the previous lemma, for every compact set $Q\subset\R^n$ and every $\delta>0$, on an event of probability $1$ we have, for all $j$ sufficiently large and every $k$ such that $\lambda_{j,k}\cap Q\neq\emptyset$,
\begin{align}\label{eqn:minorationduproduitscalaire}
|\langle f, \Psi_{p,j,k}^{\perp}\rangle|&\geq 2^{-jH}|Z_{j,k}(\delta)|- C_Q^* j 2^{-j[H+\delta (2\alpha-n)/2]}. 
\end{align}
On the other hand, $Z_{j,k}(\delta)$ is a centered Gaussian random variable with variance

\begin{align*}
V_{1,j,k}&:=\E[|Z_{j,k}(\delta)|^2]\\
&=\sum_{p'=1}^{2^n-1}\Bigg(
\sum_{|k'-k|<2^{j\delta}}|\widehat\Lambda_{p',p}(k'-k)|^2
+2^{-2H}\sum_{|k'-2k|<2^{j\delta}}|\widehat\Lambda^-_{p',p}(k'-2k)|^2\\
&\hspace{3cm}+2^{2H}\sum_{|k'-k/2|<2^{j\delta}}|\widehat\Lambda^+_{p',p}(k'-k/2)|^2\Bigg).
\end{align*}

Remark that $\sup_{j,k} V_{1,j,k} <+\infty$ as the $\widehat{\Lambda}$ functions are in $\mathscr{S}^\alpha$ and thus decay sufficiently rapidly. Moreover, it will be useful to remark that $V_{1,j,k}$ is strictly positive for $j$ sufficiently large. That is the goal of the following lemma. 

\begin{lemma} \label{lemme:variancepos}
    There exists $J_0>0$ such that 
    \[V_1 = \inf_{j\geq J_0,k\in \Z^n} V_{1,j,k}>0.\]
\end{lemma}

\begin{proof}
Remark that the functions $\widehat{\Lambda}$ appearing in
$V_{1,j,k}$ can be seen, up to the factors $2^{\pm H}$,
as the wavelet coefficients of
\[
g_k(\xi):=|\xi|^{H+n/2}\, \overline{K(\xi)}\,  \widehat{\psi}_{p,0,k}(\xi)
\]
in the basis $(\widehat{\psi}_{p',r,k'})_{p',r,k'}$.
More precisely, using binary self-similarity and Plancherel's
formula, we have, for every $k,k'\in\Z^n$ and
$p'\in\{1,\dots,2^n-1\}$,
\begin{itemize}
\item[$\bullet$] $\langle \widehat{\psi}_{p',r,k'},g_k\rangle=0$ as soon as $|r|>1$,
\item[$\bullet$] $\langle \widehat{\psi}_{p',0,k'},g_k\rangle= \widehat{\Lambda}_{p',p}(k'-k)$,
\item[$\bullet$]$\langle \widehat{\psi}_{p',1,k'},g_k\rangle =2^{-}\widehat{\Lambda}^-_{p',p}(k'-2k)$,
\item[$\bullet$] $\langle \widehat{\psi}_{p',-1,k'},g_k\rangle 2^H\widehat{\Lambda}^+_{p',p}(k'-k/2)$.
\end{itemize}
By our choice of $p$, $g_k\neq 0$ for every $k\in\Z^n$. Moreover, since $\widehat{\psi}_{p,0,k}(\xi) =e^{-ik\cdot\xi}\widehat{\psi}_p(\xi)$, we have $\|g_k\|_{L^2(\R^n)}=\|g_0\|_{L^2(\R^n)}>0$. Thus, Parseval's identity gives, for every $k\in\Z^n$,
\[
\sum_{p',r,k'}|\langle\widehat{\psi}_{p',r,k'},g_k\rangle|^2=(2\pi)^n\|g_0\|_{L^2(\R^n)}^2>0.
\]
On the other hand, since the $\widehat{\Lambda}$ functions belong to $\mathscr{S}^\alpha$, we have, uniformly in $k\in\Z^n$,
\[
0\leq(2\pi)^n\|g_0\|_{L^2(\R^n)}^2-V_{1,j,k}\leq C2^{-j\delta(2\alpha-n)}.
\]
Consequently, choosing $J_0$ sufficiently large, we obtain
\[
\inf_{j\geq J_0,\ k\in\Z^n}V_{1,j,k}
\geq\frac{(2\pi)^n}{2}\|g_0\|_{L^2(\R^n)}^2>0,
\]
which concludes the proof.
\end{proof}

In what follows, we restrict attention to $j\geq J_0$, which does not affect any of the asymptotic statements.

\begin{remark}\label{remark:realpart}
For every $(j,k)$, let $Y_{j,k}$ be either $\Re Z_{j,k}$ or $\Im Z_{j,k}$, choosing a component whose variance is maximal. Then
\[
V_{2,j,k}:=\E[|Y_{j,k}|^2]\geq \frac{1}{2}V_{1,j,k}\geq V_2:=\frac{V_1}{2}>0.
\]
The choice is deterministic, since it only depends on the covariance of $Z_{j,k}$. Moreover, if two variables $Z_{j,k}$ and $Z_{j',k'}$ are independent, then so are $Y_{j,k}$ and $Y_{j',k'}$.
\end{remark}

\begin{remark}\label{rem:independence}
    Remark that the random variables $Z_{j,k}(\delta)$ and $Z_{j',\ell}(\delta)$ are independent as soon as $|j-j'|>2$. Indeed, $Z_{j,k}(\delta)$ is composed of independent $\mathscr{N}(0,1)$ random variables of the form $\varepsilon_{\cdot, j,\cdot},\varepsilon_{\cdot, j-1,\cdot}$ and $\varepsilon_{\cdot, j+1,\cdot}$ and by taking $|j-j'|>2$, we ensure that 
    \[\{j-1,j,j+1\}\cap \{j'-1,j',j'+1\} = \emptyset.\]
\end{remark}

From now on, we fix $\delta\in(0,1)$ and will write $Z_{j,k}$ for $Z_{j,k}(\delta)$. We can now control the sequence $(Z_{j,k})_{j,k}$ according to the moduli of continuity in \eqref{eqn:modofcont}. We have the three following lemmas which are inspired by \cite[Lemmata 3.3, 3.6 and 3.8]{Ayache2019a}.

\begin{lemma}
\label{thm : asymptotic lower bound 1}
There exists an event $\Omega^*$ of probability $1$ on which, for almost every $t\in \R^n$, one has 
\[\limsup_{j\to \infty} \dfrac{|Z_{j,k_j(t)}|}{\sqrt{\log(j)}}>0.\]
\end{lemma}
\begin{proof}
Fix $t\in \R^n$ and set $u_j=\sqrt{\log j}$ for $j\geq2$. Using Remark \ref{remark:realpart}, we work with $Y_{j,k_j(t)}$, $j\in\N$. For all such a $j$, we have 

\begin{align*}
    \p(|Y_{j,k}|> \sqrt{V_2} u_j) \geq \p(|Y_{j,k}|> \sqrt{V_{2,j,k}} u_j)\geq C \dfrac{1}{\sqrt{\log(j)}} e^{-\log(j)/2} = C \dfrac{1}{\sqrt{j\log(j)}},
\end{align*}
where we used the fact $\p(|X|\geq c)\sim (2\pi^{-1})^{1/2}c^{-1}e^{-c^2/2} $ when $c\to \infty$ for any $X\sim \mathscr{N}(0,1)$.
It follows that
\[\sum_{j\in 3\N}  \p(|Y_{j,k_j(t)}|> \sqrt{V_2} u_j) = +\infty,\]
Thus using the independence (see Remark \ref{rem:independence}) and the second Borel-Cantelli Lemma, it follows that there exists infinitely many $j>0$ such that 
\[|Y_{j,k_j(t)}| \geq \sqrt{V_2} u_j>0. \]
The conclusion follows from Fubini's theorem.
\end{proof}

\begin{lemma}
\label{thm : asymptotic lower bound 2}
There exists an event $\Omega^*$ of probability $1$ on which, for every
non-empty open set $I\subset\R^n$, there exists $t\in I$ such that
\[
\limsup_{j\to\infty}\frac{|Z_{j,k_j(t)}|}{\sqrt{j}}>0.
\]
\end{lemma}

\begin{proof}
The proof follows the same lines as that of
\cite[Lemma 3.6]{Ayache2019a}. We only indicate the modifications required
in the present setting. As before, using Remark \ref{remark:realpart}, we work with the real
Gaussian variables $Y_{j,k}$. Fix $\delta<a<1$ and choose $C>0$ such that
\[
C^2<2n(a-\delta)\log(2)V_2.
\]
For $j\in\N$, set
\[
R_j=\lceil2^{j\delta}\rceil,\quad
M_j=\lfloor2^{(a-\delta)j}\rfloor
\quad \text{and} \quad
L_j=\left\lfloor\frac{2^j}{4R_jM_j}\right\rfloor.
\]
For $l\in\{0,\dots,L_j-1\}^n$, let
\[
K_{j,l} =\prod_{i=1}^n \{l_iM_j,\dots,(l_i+1)M_j-1\}
\]
and define
\[
\mathcal F_{j,l} = \left\{\max_{k\in K_{j,l}}|Y_{j,4R_jk}|\ge C\sqrt j\right\}.
\]
The random variables appearing in this maximum are independent. Indeed, $Z_{j,k}$ only involves Gaussian coefficients whose spatial indices are at distance less than $2^{j\delta}$ from $k$, after rescaling
to the central scale $j$. Since two distinct indices of the form $4R_jk$ are separated by at least $4R_j\ge4\,2^{j\delta}$ in at least one coordinate, the corresponding families of Gaussian coefficients
are disjoint. Let
\[
q_j=\p\left(|\varepsilon|\ge\frac{C}{\sqrt{V_2}}\sqrt j\right),\quad \text{ for }
\varepsilon\sim\mathcal N(0,1).
\]
Since $V_{2,j,k}\ge V_2$, we have
\[
\p\left(|Y_{j,k}|<C\sqrt j\right)\le 1-q_j.
\]
Moreover, the classical Gaussian tail estimate gives, for $j$ large enough,
\[
q_j
\ge
c\,j^{-1/2}
\exp\left(-\frac{C^2}{2V_2}j\right).
\]
Therefore, using the independence above,
\[
\p(\mathcal F_{j,l}^c)
\le
(1-q_j)^{M_j^n}
\le
\exp(-M_j^nq_j).
\]
Our choice of $C$  and of $M_j$ implies that there exists $\eta>0$ such that
\[
M_j^nq_j
\ge
c\,j^{-1/2}e^{\eta j}
\]
for all sufficiently large $j$. As $L_j^n\le C'2^{n(1-a)j}$, we obtain 
\[
\p\left(
\bigcup_{l\in\{0,\dots,L_j-1\}^n}
\mathcal F_{j,l}^c
\right)
\le
C'2^{n(1-a)j}
\exp\left(-c\,j^{-1/2}e^{\eta j}\right).
\]
The right-hand side is summable in $j$. Hence, by the Borel-Cantelli
Lemma, almost surely, for all sufficiently large $j$, every block
$K_{j,l}$ contains an index $k$ such that
\[
|Y_{j,4R_jk}|\ge C\sqrt j.
\]
For $\omega$ in this event, set
\[
G_j(\omega)
=
\left\{
k\in\{0,\dots,2^j-1\}^n:
|Y_{j,k}(\omega)|\ge C\sqrt j
\right\}
\]
and
\[
U_j(\omega)
=
\bigcup_{k\in G_j(\omega)}
\operatorname{int}(\lambda_{j,k}),
\qquad
O_m(\omega)=\bigcup_{j\ge m}U_j(\omega).
\]

The preceding block construction implies that $U_j(\omega)$ is
$C_0 2^{-j(1-a)}$-dense in $(0,1)^n$, for some deterministic constant
$C_0>0$. Indeed, the diameter in the original dyadic lattice of each
block of selected indices is of order $R_jM_j$, which behaves as $2^{aj}$. Since $a<1$, it follows that every $O_m(\omega)$ is open and dense in $(0,1)^n$. By Baire's Theorem,
\[
\bigcap_{m\ge0}O_m(\omega)
\]
is non-empty. If $t$ belongs to this intersection, then for arbitrarily large $j$,
\[
|Y_{j,k_j(t)}|\ge C\sqrt j.
\]
Since $|Z_{j,k}|\ge|Y_{j,k}|$, we obtain
\[
\limsup_{j\to\infty} \frac{|Z_{j,k_j(t)}|}{\sqrt j}
\ge C>0.
\]
The same argument applies to every dyadic cube. Taking the intersection of the corresponding events over the countable family of dyadic cubes gives a common event of probability $1$. Since every non-empty open subset of $\R^n$ contains a dyadic cube, the conclusion follows.
\end{proof}

 \begin{lemma}
 \label{thm : asymptotic lower bound 3}
 There exists an event $\Omega^*$ of probability $1$ on which, for every
 $t\in\R^n$,
 \[
 \limsup_{j\to\infty}|Z_{j,k_j(t)}|>0.
 \]
 \end{lemma}

\begin{proof}
Choose $c_0>0$ so small that
\[
p_0:=\sup_{j\geq J_0,\,k\in\Z^n}\p\left(|Z_{j,k}|\leq c_0\right)<2^{-3n}.
\]
Such a choice is possible since, by Remark \ref{remark:realpart},
\[
\p(|Z_{j,k}|\leq c_0)\leq \p(|Y_{j,k}|\leq c_0)
\leq \p(|N(0,V_2)|\leq c_0),
\]
uniformly in $j\geq J_0$ and $k\in\Z^n$. Fix $(J,k)\in\N\times\Z^n$, with $J\geq J_0$, and $m\in\N$.
For every dyadic cube $S$ of scale $J+3m$ contained in $\lambda_{J,k}$,
let
\[
\lambda_{J+3\ell,k_{J+3\ell}(S)},
\qquad 0\leq\ell\leq m,
\]
denote the unique cubes of the ancestral chain of $S$ at scales
$J,J+3,\dots,J+3m$. Set
\[
B_{J,k,m}(S)
=
\prod_{\ell=0}^m\mathds{1}_{\{|Z_{J+3\ell,k_{J+3\ell}(S)}|\leq c_0\}}.
\]
By Remark \ref{rem:independence}, the random variables appearing in this
product are independent. Hence
\[
\E[B_{J,k,m}(S)]\leq p_0^{m+1}.
\]
Let
\[
G_{J,k,m}
=\!\!\!\!\!\!
\sum_{\substack{S\subset\lambda_{J,k}\\
\operatorname{scale}(S)=J+3m
}}
\!\!\!\!\!\! B_{J,k,m}(S).
\]
Since there are $2^{3nm}$ such cubes,
\[
\E[G_{J,k,m}]
\leq
2^{3nm}p_0^{m+1}.
\]
Since $p_0<2^{-3n}$, the right-hand side tends to $0$ as $m\to\infty$.
Therefore, by Fatou's Lemma,
\[
0
\leq
\E\left[\liminf_{m\to\infty}G_{J,k,m}\right]\leq\liminf_{m\to\infty}\E[G_{J,k,m}]=0.
\]
Consequently,
\[
\liminf_{m\to\infty}G_{J,k,m}=0\qquad\text{almost surely}.
\]
Since $G_{J,k,m}$ is integer-valued, this implies that
\[
G_{J,k,m}=0
\]
for infinitely many $m$. Taking the intersection over the countable family
$(J,k)\in\N\times\Z^n$, $J\geq J_0$, we obtain a common event of
probability $1$ on which this property holds simultaneously for every
$(J,k)$.

Let now $t\in\R^n$ and fix arbitrarily $J\geq J_0$. Apply the previous
conclusion with $k=k_J(t)$. For infinitely many $m$, we have
$G_{J,k_J(t),m}=0$. In particular, taking
\[
S=\lambda_{J+3m,k_{J+3m}(t)},
\]
we have $B_{J,k_J(t),m}(S)=0$. Hence, for each such $m$, there exists
$\ell\in\{0,\dots,m\}$ such that
\[
|Z_{J+3\ell,k_{J+3\ell}(t)}|>c_0.
\]
Since $J$ can be chosen arbitrarily large, such indices occur at
arbitrarily large scales. Therefore,
\[
\limsup_{j\to\infty}|Z_{j,k_j(t)}|\geq c_0>0,
\]
which concludes the proof.
\end{proof}

We can finally prove that the three limits in the Theorem \ref{theoreme principal} are strictly positive.
\begin{proof}[Proof of Theorem \ref{theoreme principal}]
For $m\geq1$, set $Q_m=[-m,m]^n$. Taking the intersection of the events supplied by the estimate \eqref{eqn:minorationduproduitscalaire} for the countable family $(Q_m)_{m\geq1}$, and intersecting it with the common events from Proposition \ref{prop : maj produit scalaire} and the three previous lemmas, we obtain an event $\widehat\Omega$ of probability $1$ on which all these estimates can be used simultaneously.

Fix $t\in\R^n$ and choose $m$ such that $t\in Q_m$. Then, for every sufficiently large $j$, $\lambda_{j,k_j(t)}\cap Q_m\neq\emptyset$, and therefore
\begin{align*}
2^{-jH}|Z_{j,k_j(t)}|-C_{Q_m}^*j2^{-j[H+\delta(2\alpha-n)/2]}
&\leq |\langle f,\Psi^{\perp}_{p,j,k_j(t)}\rangle|\\
&\leq C'\left(\sup\left\{\frac{|f(t)-f(s)|}{\sigma(|t-s|)}:0<|t-s|<c2^{-j/2}\right\}\sigma(2^{-j})+2^{-j\gamma}\right),
\end{align*}
where $\gamma>H$. For each of the three moduli in \eqref{eqn:modofcont},
\[
\frac{j2^{-j[H+\delta(2\alpha-n)/2]}}{\sigma_i(2^{-j})}\longrightarrow0,
\qquad
\frac{2^{-j\gamma}}{\sigma_i(2^{-j})}\longrightarrow0.
\]
Consequently, Lemmas \ref{thm : asymptotic lower bound 1}, \ref{thm : asymptotic lower bound 2} and \ref{thm : asymptotic lower bound 3} yield respectively
    \begin{itemize}
        \item[$\bullet$] for almost every $t\in \R^n$, almost surely,
        \[0< \limsup_{j \to +\infty} \dfrac{ 2^{-jH} |Z_{j,k_j(t)}|}{\sigma_1(2^{-j})} \leq C \lim_{j \to +\infty} \sup\left\{\dfrac{|f(t)-f(s)|}{\sigma_1(|t-s|)}:0<|t-s|<c2^{-j/2}\right\}; \]
        \item[$\bullet$] Almost surely, for all non-empty open set $I\subset \R^n$, there exists $t \in I$ such that 
        \[0< \limsup_{j \to +\infty} \dfrac{ 2^{-jH} |Z_{j,k_j(t)}|}{\sigma_2(2^{-j})} \leq C \lim_{j \to +\infty} \sup\left\{\dfrac{|f(t)-f(s)|}{\sigma_2(|t-s|)}:0<|t-s|<c2^{-j/2}\right\}; \]
        \item[$\bullet$] Almost surely, for every $t\in \R^n$,
        \[0< \limsup_{j \to +\infty} \dfrac{ 2^{-jH} |Z_{j,k_j(t)}|}{\sigma_3(2^{-j})} \leq C \lim_{j \to +\infty} \sup\left\{\dfrac{|f(t)-f(s)|}{\sigma_3(|t-s|)}:0<|t-s|<c2^{-j/2}\right\}. \]
        \end{itemize}

    Now, using similar arguments as in \cite[Propositions 3.2, 3.5,3.8] {Esser2022}, we already know that 
        \begin{itemize}
        \item[$\bullet$] for almost every $t\in \R^n$, almost surely,
        \[ \lim_{j \to +\infty} \sup\left\{\dfrac{|f(t)-f(s)|}{\sigma_1(|t-s|)}:0<|t-s|<c2^{-j/2}\right\}<+\infty; \]
        \item[$\bullet$] Almost surely, for every $t \in \R^n$, 
        \[ \lim_{j \to +\infty} \sup\left\{\dfrac{|f(t)-f(s)|}{\sigma_2(|t-s|)}:0<|t-s|<c2^{-j/2}\right\}< +\infty; \]
        \item[$\bullet$] Almost surely, for all non-empty open set $I\subset \R^n$, there exists $t \in I$ such that
        \[\lim_{j \to +\infty} \sup\left\{\dfrac{|f(t)-f(s)|}{\sigma_3(|t-s|)}:0<|t-s|<c2^{-j/2}\right\}< + \infty. \]
        \end{itemize}
Since $X-f$ is differentiable and $H<1$, its increments are negligible with respect to the three moduli $\sigma_1,\sigma_2,\sigma_3$. Therefore the corresponding limsup statements for $f$ are equivalent to those for $X$. Now, combining the corresponding bullets points gives the result respectively for the ordinary, fast and slow points. More precisely, for the rapid points, remark that for the lower bound, we only know that there exists $t$ while the upper bound works for every $t$ and thus, the point $t$ given by the lower bound is a rapid point. For the slow points, the same reasoning applies by interchanging the upper and lower bound while for the ordinary points, the conclusion follows as the intersection of two almost sure sets is also almost sure.
\end{proof}

\section{Proof of Theorem \ref{thm:moving}}

Thanks to the development in Section \ref{sec:proofmain1}, one can give a moving average representation of the field $X$.
Indeed, from the proof of Proposition \ref{prop wavelet development via dist}, we know that 
\begin{equation}
    \label{eq9}
(e^{it\cdot\xi}-1) K(\xi) = (2\pi)^{-n/2}\sum_{p,j,k} 2^{-jH}\Psi_p(-k)(e^{it\cdot\xi}-1) \overline{\widehat{\psi}}_{p,j,k}(\xi)
\end{equation}
where the equality is understood in $L^2(\R^n)$; the expansion of $K$ itself is locally uniform away from the origin. That is the main ingredient to prove Theorem \ref{thm:moving}

\begin{proof}[Proof of Theorem \ref{thm:moving}]
    From  the identity \eqref{eqn:intwithfouriermeasure} and the representation \eqref{eq9}, we directly obtain that, for any $t \in \R^n$, 
    \begin{align*}
    X_t &= \int_{\R^n} \sum_{p,j,k}2^{-jH} \Psi_p(-k) \left( \overline{\psi}_{p,j,k}(u-t) -\overline{\psi}_{p,j,k}(u)\right) \, dW(u)\\
    &= \int_{\R^n}\varphi(u-t)-\varphi(u)\, dW(u),
    \end{align*}
which is well defined as $u\mapsto \varphi(u-t)-\varphi(u)\in L^2(\R^n)$ as it is the Fourier transform of an $L^2(\R^n)$ function. More precisely, the $L^2$ convergence of the series of differences was established in Proposition \ref{prop wavelet development via dist}; the locally uniform convergence of the series defining $\varphi$ away from $0$ identifies its $L^2$ limit with $\varphi(\cdot-t)-\varphi$ almost everywhere. Thus we only need to show that $\varphi$ is well-defined. 
Let us fix $j \in \Z$, $p\in \{1,\dots, 2^n-1\}$ and $\xi \in \R^n\setminus\{0\}$, we have
\begin{align*}
    \sum_{k\in \Z^n} 2^{-jH} |\Psi_p(-k)| |\psi_{p,j,k}(\xi)| &= \sum_{|k|\leq 2^{j-1}|\xi|} 2^{-jH} |\Psi_p(-k)| |\psi_{p,j,k}(\xi)|+\sum_{|k|>2^{j-1}|\xi|} 2^{-jH}|\Psi_p(-k)| |\psi_{p,j,k}(\xi)|.
\end{align*}
Let us bound the first sum. Using Lemma \ref{lemma:conv_serie}, we have, for any $N>0$,
\begin{align*}
   \sum_{|k|\leq 2^{j-1}|\xi|} 2^{-jH} |\Psi_p(-k)| |\psi_{p,j,k}(\xi)|&\leq \sum_{|k|\leq 2^{j-1}|\xi|} 2^{-jH}2^{jn/2} |\Psi_p(-k)| \dfrac{C}{(1+|2^j\xi-k|)^N} \\
  & \leq \sum_{|k|\leq 2^{j-1}|\xi|} 2^{-jH}2^{jn/2} |\Psi_p(-k)| \dfrac{C}{(1+2^j|\xi|-|k|)^N} \\
   &\leq \sum_{|k|\leq 2^{j-1}|\xi|} 2^{-jH}2^{jn/2} |\Psi_p(-k)| \dfrac{C}{(1+2^{j-1}|\xi|)^N} \\
   &\leq C' 2^{-jH}2^{jn/2} \dfrac{C}{(1+2^{j-1}|\xi|)^N}
\end{align*}
as $
\sum_{k \in \Z^n} |\Psi_p(-k)| < + \infty.
$

On the other side, for the second sum, using again Lemma \ref{lemma:conv_serie}, we have, for any $n<N\leq \alpha$
\begin{align*}
    \sum_{|k|>2^{j-1}|\xi|} 2^{-jH}|\Psi_p(-k)| |\psi_{p,j,k}(\xi)|&\leq \sum_{|k|>2^{j-1}|\xi|} 2^{-jH} \dfrac{C}{(1+|k|)^N} |\psi_{p,j,k}(\xi)|\\
    &\leq \sum_{|k|>2^{j-1}|\xi|} 2^{-jH} \dfrac{C}{(1+2^{j-1}|\xi|)^N} |\psi_{p,j,k}(\xi)|\\
    &\leq  C' 2^{-jH}2^{jn/2}\dfrac{C}{(1+2^{j-1}|\xi|)^N}.
\end{align*}

To conclude the proof, it suffices to show that, for $n<N\leq \alpha$, the series 
\[\sum_{j\in \Z} 2^{-jH}2^{jn/2}\dfrac{1}{(1+2^{j-1}|\xi|)^N} \]
converges uniformly on every compact of $\R^n\setminus\{0\}$. Fix $F$ such a compact set. Then $\inf |F| >0$ and for every $\xi \in F$, one has
$$
\dfrac{1}{(1+2^{j-1}|\xi|)^N} \leq \dfrac{1}{(2^{j-1}|\xi|)^N} \leq \frac{1}{(2^{j-1}\inf |F|)^N}.
$$
It follows that for the terms associated with  positive values of $j$, we have
\[\sup_{\xi \in F} \sum_{j\geq0} 2^{-jH}2^{jn/2}\dfrac{1}{(1+2^{j-1}|\xi|)^N}\leq (\inf |F|)^{-N} \sum_{j\geq 0} 2^{-jH}2^{jn/2}2^{(1-j)N} <\infty.\]
For the negative $j$, using the assumption $H<n/2$, one obtains 
\[\sup_{\xi \in F} \sum_{j<0} 2^{-jH}2^{jn/2}\dfrac{1}{(1+2^{j-1}|\xi|)^N}\leq  \sum_{j<0} 2^{-jH}2^{jn/2}<\infty.\]
\end{proof}

\noindent\textbf{Acknowledgements.} 
This work was supported by grant No. 1.A.679.26F. from the Fonds de la Recherche
Scientifique (F.R.S.-FNRS, Belgium) awarded to T. Devos.
C. Esser acknowledges financial support from the Fonds de la
Recherche Scientifique (F.R.S.-FNRS, Belgium) through CFR funding
No. J.0067.26, and from the Fonds Thelam,
managed by the King Baudouin Foundation (Belgium).
L. Loosveldt acknowledges financial support from the Fonds de la
Recherche Scientifique (F.R.S.-FNRS, Belgium) through CFR funding No. J.0136.26 ``NOISE''.

\medskip
\noindent\textbf{Statement on AI use.} The authors used ChatGPT (OpenAI, GPT-5.6 Sol and GPT-6.2 Astra) in the preparation of this manuscript, for language editing, improving the clarity of the exposition, and clarifying certain mathematical arguments. All AI-assisted content was carefully reviewed and, where necessary, revised by the authors, who take full responsibility for the final content of the manuscript.

\bibliography{biblio}{}

\begin{thebibliography}{10}

\bibitem{ABRY20062326}
P.~Abry and V.~Pipiras.
\newblock Wavelet-based synthesis of the rosenblatt process.
\newblock {\em Signal Processing}, 86(9):2326--2339, 2006.
\newblock Special Section: Signal Processing in UWB Communications.

\bibitem{abrysellan96}
P.~Abry and F.~Sellan.
\newblock The wavelet-based synthesis for fractional {B}rownian motion proposed
  by {F}. {S}ellan and {Y}. {M}eyer: remarks and fast implementation.
\newblock {\em Appl. Comput. Harmon. Anal.}, 3(4):377--383, 1996.

\bibitem{Ayache2019}
A.~Ayache.
\newblock {\em Multifractional stochastic fields}.
\newblock World Scientific Publishing Co. Pte. Ltd., Hackensack, NJ.
\newblock Wavelet strategies in multifractional frameworks.

\bibitem{ayache2002}
A.~Ayache.
\newblock The generalized multifractional field: a nice tool for the study of
  the generalized multifractional {B}rownian motion.
\newblock {\em J. Fourier Anal. Appl.}, 8(6):581--601, 2002.

\bibitem{Ayache2019a}
A.~Ayache, C.~Esser, and T.~Kleyntssens.
\newblock Different possible behaviors of wavelet leaders of the {B}rownian
  motion.
\newblock 150:54--60.

\bibitem{ayachehamonier}
A.~Ayache and J.~Hamonier.
\newblock Wavelet series representation for multifractional multistable
  {R}iemann-{L}iouville process.
\newblock {\em Math. Methods Econ. Finance}, 13/14:1--17, 2018/19.

\bibitem{ayachehamonierloosveldt}
A.~Ayache, J.~Hamonier, and L.~Loosveldt.
\newblock Wavelet-type expansion of generalized {H}ermite processes with rate
  of convergence.
\newblock {\em Constr. Approx.}, 61(3):535--601, 2025.

\bibitem{ayachehamonierloosveldt2}
A.~Ayache, J.~Hamonier, and L.~Loosveldt.
\newblock Numerical simulation of generalized {H}ermite processes.
\newblock {\em Ann. Appl. Probab.}, 36(4):3601--3634, 2026.

\bibitem{ayachelegerpontier}
A.~Ayache, S.~Leger, and M.~Pontier.
\newblock Les ondelettes \`a{} la conqu\^ete du drap brownien fractionnaire.
\newblock {\em C. R. Math. Acad. Sci. Paris}, 335(12):1063--1068, 2002.

\bibitem{ayachelouckx}
A.~Ayache and C.~Louckx.
\newblock Harmonizable multifractional stable field: sharp results on sample
  path behavior.
\newblock {\em Stochastic Process. Appl.}, 186:Paper No. 104638, 35, 2025.

\bibitem{ayacheroueffxiao}
A.~Ayache, F.~Roueff, and Y.~Xiao.
\newblock Linear fractional stable sheets: wavelet expansion and sample path
  properties.
\newblock {\em Stochastic Process. Appl.}, 119(4):1168--1197, 2009.

\bibitem{ayacheshiehxiao}
A.~Ayache, N.~Shieh, and Y.~Xiao.
\newblock Wavelet series representation and geometric properties of
  harmonizable fractional stable sheets.
\newblock {\em Stochastics}, 92(1):1--23, 2020.

\bibitem{ayachetaqqu2003}
A.~Ayache and M.~S. Taqqu.
\newblock Rate optimality of wavelet series approximations of fractional
  {B}rownian motion.
\newblock {\em J. Fourier Anal. Appl.}, 9(5):451--471, 2003.

\bibitem{ayachetaqqu2005}
A.~Ayache and M.~S. Taqqu.
\newblock Multifractional processes with random exponent.
\newblock {\em Publ. Mat.}, 49(2):459--486, 2005.

\bibitem{Ayachexiao2026}
A.~Ayache and Y.~Xiao.
\newblock On directional sample path behavior of harmonizable fractional stable
  sheets.
\newblock {\em Journal of Mathematical Analysis and Applications},
  562(2):130681, 2026.

\bibitem{Baek2014}
C.~Baek, G.~Didier, and V.~Pipiras.
\newblock On integral representations of operator fractional {B}rownian fields.
\newblock 92:190--198.

\bibitem{Bierme2026}
H.~Bierm\'e, P.~Carr\'e, C.~Lacaux, and C.~Launay.
\newblock Gaussian random fields and monogenic images.
\newblock 80:Paper No. 101814, 27.

\bibitem{Bierme2007}
H.~Bierm\'e, M.~M. Meerschaert, and H.~Scheffler.
\newblock Operator scaling stable random fields.
\newblock 117(3):312--332.

\bibitem{Bierme2015}
H.~Bierm\'e, L.~Moisan, and F.~Richard.
\newblock A turning-band method for the simulation of anisotropic fractional
  {B}rownian fields.
\newblock {\em J. Comput. Graph. Statist.}, 24(3):885--904, 2015.

\bibitem{Bierme2009}
H.~Bierm\'e, F.~Richard, M.~Rachidi, and C.~Benhamou.
\newblock Anisotropic texture modeling and applications to medical image
  analysis.
\newblock In {\em Mathematical methods for imaging and inverse problems},
  volume~26 of {\em ESAIM Proc.}, pages 100--122. EDP Sci., Les Ulis, 2009.

\bibitem{Bonami2003}
A.~Bonami and A.~Estrade.
\newblock Anisotropic analysis of some {G}aussian models.
\newblock {\em J. Fourier Anal. Appl.}, 9(3):215--236, 2003.

\bibitem{Davy2023}
L.~Davy, N.~Pustelnik, and P.~Abry.
\newblock Combining dual-tree wavelet analysis and proximal optimization for
  anisotropic scale-free texture segmentation.
\newblock In {\em ICASSP 2023 - 2023 IEEE International Conference on
  Acoustics, Speech and Signal Processing (ICASSP)}, pages 1--5, 2023.

\bibitem{Daw2022}
L.~Daw and L.~Loosveldt.
\newblock Wavelet methods to study the pointwise regularity of the generalized
  {R}osenblatt process.
\newblock 27:Paper No. 152, 45.

\bibitem{Escande2024}
P.~Escande and F.~J.~P. Richard.
\newblock Full inference for the anisotropic fractional {B}rownian field.
\newblock {\em Theory Probab. Math. Statist.}, (110):13--29, 2024.

\bibitem{Esser2024}
C.~Esser and L.~Loosveldt.
\newblock On the pointwise regularity of the multifractional {B}rownian motion
  and some extensions.
\newblock (110):55--73.

\bibitem{Esser2022}
C.~Esser and L.~Loosveldt.
\newblock Slow, ordinary and rapid points for {G}aussian wavelets series and
  application to fractional {B}rownian motions.
\newblock 19(2):1471--1495.

\bibitem{herbin}
E.~Herbin.
\newblock From {$N$} parameter fractional {B}rownian motions to {$N$} parameter
  multifractional {B}rownian motions.
\newblock {\em Rocky Mountain J. Math.}, 36(4):1249--1284, 2006.

\bibitem{Janson1997}
S.~Janson.
\newblock {\em Gaussian Hilbert Spaces}.
\newblock Cambridge Tracts in Mathematics. Cambridge University Press, 1997.

\bibitem{kahane}
J.~Kahane.
\newblock {\em Some random series of functions}, volume~5 of {\em Cambridge
  Studies in Advanced Mathematics}.
\newblock Cambridge University Press, Cambridge, second edition, 1985.

\bibitem{lemariemeyer}
P.~G. Lemari\'e and Y.~Meyer.
\newblock Ondelettes et bases hilbertiennes.
\newblock {\em Rev. Mat. Iberoamericana}, 2(1-2):1--18, 1986.

\bibitem{Maejima1999}
M.~Maejima and K.~Sato.
\newblock Semi-selfsimilar processes.
\newblock {\em J. Theoret. Probab.}, 12(2):347--373, 1999.

\bibitem{meyerselaantaqqu}
Y.~Meyer, F.~Sellan, and M.~S. Taqqu.
\newblock Wavelets, generalized white noise and fractional integration: the
  synthesis of fractional {B}rownian motion.
\newblock {\em J. Fourier Anal. Appl.}, 5(5):465--494, 1999.

\bibitem{pipiras05}
V.~Pipiras.
\newblock Wavelet-based simulation of fractional {B}rownian motion revisited.
\newblock {\em Appl. Comput. Harmon. Anal.}, 19(1):49--60, 2005.

\bibitem{Polisano2015}
K.~Polisano, M.~Clausel, V.~Perrier, and L.~Condat.
\newblock Mod{\'e}lisations de textures par champ gaussien {\`a} orientation
  locale prescrite.
\newblock Preprint, {arXiv}:1503.06716 [math.{CA}] (2015), 2015.

\bibitem{Xiao2021TenLectures}
Y.~Xiao.
\newblock Ten lectures on ``gaussian random fields, spdes, fractals, and
  extremes''.
\newblock Lecture notes, Michigan State University, 2021.

\end{thebibliography}
\bibliographystyle{plain}

\end{document}